\documentclass[11pt,a4paper]{article}
\usepackage{amsfonts,amsmath,amssymb,accents,amsthm}%,mathtools}
\usepackage{verbatim}
\usepackage{hyperref}
\usepackage[normalem]{ulem}%for \sout,\xout

\usepackage{tikz}
\usetikzlibrary{calc,arrows}

\newcommand{\ii}{\imath}
\newcommand{\jj}{\jmath}
\newcommand{\dott}{\bullet}
\newcommand{\both}{\sharp}

\newcommand{\status}{}

\newcommand{\file}{}

\newcommand{\detail}[1]{\par\noi{\bf [Proof detail\ }{#1}
\hfill{\bf ]}\par\noi\hspace{-4pt}}
\renewcommand{\detail}[1]{}

\newcommand{\dis}{\displaystyle}
\newcommand{\txt}{\textstyle}

\newcommand{\noi}{\noindent}
\newcommand{\halmos}{\rule{1ex}{1.4ex}}
\newcommand{\QED}{\nopagebreak{\hspace*{\fill}$\halmos$\medskip}}

\newcommand{\med}{\medskip}
\newcommand{\quand}{\quad\mbox{and}\quad}

\newtheoremstyle{mythm}% name
  {}%      Space above
  {}%      Space below
  {\itshape}%         Body font
  {}%         Indent amount (empty = no indent, \parindent = para indent)
  {\bfseries}% Thm head font
  {}%        Punctuation after thm head
  {.5em}%     Space after thm head: " " = normal interword space;
  {#1 #2 \thmnote{(#3)}}%  Thm head spec (can be left empty, meaning `normal')

\theoremstyle{mythm}
\newtheorem{theorem}{Theorem}%[section]
\newtheorem{proposition}[theorem]{Proposition}
\newtheorem{lemma}[theorem]{Lemma}
\newtheorem{exercise}[theorem]{Exercise}
\newtheorem{corollary}[theorem]{Corollary}
\newtheorem{conjecture}[theorem]{Conjecture}

\newtheorem{counterex}[theorem]{Counterexample}

\newcommand{\bt}{\begin{theorem}}
\newcommand{\et}{\end{theorem}}
\newcommand{\bl}{\begin{lemma}}
\newcommand{\el}{\end{lemma}}
\newcommand{\bp}{\begin{proposition}}
\newcommand{\ep}{\end{proposition}}
\newcommand{\bcor}{\begin{corollary}}
\newcommand{\ecor}{\end{corollary}}
\newcommand{\br}{\begin{remark}\rm}
\newcommand{\er}{\end{remark}}
\newcommand{\bcon}{\begin{conjecture}}
\newcommand{\econ}{\end{conjecture}}
\newcommand{\bex}{\begin{exercise}}
\newcommand{\eex}{\end{exercise}}
\newcommand{\bcou}{\begin{counterex}}
\newcommand{\ecou}{\end{counterex}}

\newenvironment{Proof}[1][]{\noi\textbf{Proof #1}}{\QED}
\newcommand{\bpro}{\begin{Proof}}
\newcommand{\epro}{\end{Proof}}

\newcommand{\be}{\begin{equation}}
\newcommand{\ee}{\end{equation}}
\newcommand{\ba}{\begin{array}}
\newcommand{\ea}{\end{array}}
\newcommand{\bc}{\be\begin{array}{r@{\,}c@{\,}l}}
\newcommand{\ec}{\end{array}\ee}

\newcommand{\bet}{\beta}
\newcommand{\ga}{\gamma}
\newcommand{\Ga}{\Gamma}
\newcommand{\de}{\delta}
\newcommand{\De}{\Delta}
\newcommand{\eps}{\varepsilon}

\newcommand{\La}{\Lambda}
\newcommand{\sig}{\sigma}

\newcommand{\Om}{\Omega}

\newcommand{\Ai}{{\cal A}}

\newcommand{\Di}{{\cal D}}

\newcommand{\Fi}{{\cal F}}
\newcommand{\Gi}{{\cal G}}

\newcommand{\Ki}{{\cal K}}

\newcommand{\Pc}{{\cal P}}

\newcommand{\Ri}{{\cal R}}

\newcommand{\R}{{\mathbb R}}

\newcommand{\Z}{{\mathbb Z}}

\newcommand{\E}{{\mathbb E}}
\renewcommand{\P}{{\mathbb P}}

\newcommand{\li}{\langle}
\newcommand{\re}{\rangle}

\newcommand{\desd}{\ensuremath{\Leftrightarrow}}
\newcommand{\volgt}{\ensuremath{\Rightarrow}}
\newcommand{\up}{\uparrow}
\newcommand{\down}{\downarrow}
\newcommand{\sub}{\subset}
\newcommand{\beh}{\backslash}

\newcommand{\ti}{\tilde}
\newcommand{\dgg}{\dagger}
\newcommand{\ov}{\overline}

\newcommand{\subb}[2]{_{\ba{c}\scriptstyle{#1}\\[-.15cm]\scriptstyle{#2}\ea}}

\newcommand{\di}{\mathrm{d}}
\newcommand{\half}{{[0,\infty)}}
\newcommand{\expo}{\mbox{\large\it e}}
\newcommand{\ex}[1]{\expo^{\,\textstyle{#1}}}

\newcommand{\Xb}{{\mathbf X}}
\newcommand{\Yb}{{\mathbf Y}}

\newcommand{\block}{\tikz[baseline=-3pt]{
  \draw[line width=3pt] (0,0)++(-0.1,0)--++(0.2,0);}\hspace{4pt}}

\begin{document}

%numbering formulas within sections
\makeatletter\@addtoreset{equation}{section}
\makeatother\def\theequation{\thesection.\arabic{equation}} 

%alternative layout for enumerate lists.
\renewcommand{\labelenumi}{{\rm (\roman{enumi})}}

\title{Pathwise duals of monotone and additive Markov processes}
\author{
Anja~Sturm
\footnote{Institute for Mathematical Stochastics,
Georg-August-Universit\"at G\"ottingen,
Goldschmidtstr.~7,
37077 G\"ottingen,
Germany;
asturm@math.uni-goettingen.de}
\and
Jan~M.~Swart
\footnote{Institute of Information Theory and,
Automation of the ASCR (\' UTIA),
Pod vod\'arenskou v\v e\v z\' i 4,
18208 Praha 8,
Czech Republic;
swart@utia.cas.cz}
}

\date{{\small\file}\today}

\maketitle%\vspace{-.8cm}
\status

\begin{abstract}\noi
This paper develops a systematic treatment of monotonicity-based pathwise
dualities for Markov processes taking values in partially ordered sets. We
show that every Markov process that takes values in a finite partially ordered
set and whose generator can be represented in monotone maps has a pathwise
dual process. In the special setting of attractive spin systems this has been
discovered earlier by Gray.

We show that the dual simplifies a lot when the state space is a
lattice (in the order-theoretic meaning of the word) and all monotone
maps satisfy an additivity condition. This leads to a unified treatment of
several well-known dualities, including Siegmund's dual for processes with a
totally ordered state space, duality of additive spin systems, and a duality
due to Krone for the two-stage contact process, and allows for the
construction of new dualities as well.

We show that the well-known representation of additive spin systems in terms
of open paths in a graphical representation can be generalized to additive
Markov processes taking values in general lattices, but for the process and
its dual to be representable on the same underlying space, we need to assume
that the lattice is distributive.

In the final section, we show how our results can be generalized from finite
state spaces to interacting particle systems with finite local state spaces.
\end{abstract}
\vspace{.5cm}

\noi
{\it Acknowledgement.} Work sponsored by GA\v{C}R grant P201/12/2613 and 
by DFG priority programme 1590 grant STU527/1-1.

%82C22   	Interacting particle systems
%
%60J27   	Continuous-time Markov processes on discrete state spaces
%06A06   	Partial order, general
%06B10   	Ideals, congruence relations (06BXX = Lattices)

%\newpage

%{\setlength{\parskip}{-2pt}\tableofcontents}

\newpage
\setcounter{tocdepth}{2}
\tableofcontents

\newpage
\section{Introduction}\label{S:intro}

In the present paper, we systematically investigate the link between Markov
process duality and order theory. Here, we focus on the concept of ``pathwise'' 
duality which is stronger than the more often considered ``distributional'' concept of duality.

The general problem of systematically finding Markov processes $X$ and $Y$
with state spaces $S$ and $T$ that are dual to each other with respect to a
duality function has interested a large number of researchers. In
\cite{SL95,SL97,Sud00}, a systematic treatment is given of dualities for
nearest-neighbor interacting particle systems where the duality function is of
a special ``local'' form. In \cite{GKRV09,CGGF15}, dual Markov processes are
linked to different representations of the same Lie algebra. In \cite{Moe13},
the existence of duality functions and duals is related to cones of the Markov
process. The paper \cite{JK14} investigates different notions of duality,
distinguishing in particular the ``distributional'' concept, which is the
usual modern definition of Markov process duality, and the ``pathwise''
concept.

Our unified treatment includes some of the oldest known forms of Markov
process duality. Building on earlier work, Siegmund \cite{Sie76} proved that
(almost) every monotone Markov process taking values in a totally ordered set
has a dual, which is also a monotone Markov process, taking values in (almost)
the same set. Here, the duality function is given by the indicator function on
the event that the two dual processes are ordered.

Examples of such dualities (for example between Brownian motion with
reflection and absorption at the origin) were known at least 28 years earlier;
see \cite{JK14} and also \cite[Section~II.3]{Lig85} or \cite[Section~5]{DF90}
for a (short) historical overview.

Moving away from totally ordered spaces, similar ideas have been considered on
partially ordered spaces. For example, a systematic treatment of a
generalization of Siegmund's duality to processes taking values in $\R^d$
equipped with the coordinatewise (Pareto) order has been given in
\cite{KL13}. Also, spin systems are Markov processes taking values in the set
of all subsets of a countable underlying space. Such a set of subsets,
equipped with the order of set inclusion, is a partially ordered set. It was
known since the early 1970ies that certain spin systems such as the contact
and voter models have duals. Harris \cite{Har78} showed that the essential
feature of these models is that they are additive and that there is a
percolation picture behind each additive system. The idea was further
formalized in Griffeath's monograph \cite{Gri79} where the term ``percolation
substructure'' was coined. Some years later, Gray \cite{Gra86} introduced a
more general, but also more complicated duality for spin systems that are
monotone but not necessarily additive, which includes the dualities of
additive systems as a special case.

The main contribution of the present paper is that we replace the specific
examples of partially ordered sets mentioned above (totally ordered sets
respectively the set of all subsets of another set) by a general partially
ordered set and that we focus on pathwise duality. For technical simplicity,
we will mostly concentrate on finite spaces, but we also show how the theory
can be generalized to infinite products of finite partially ordered sets,
which allows our results to be applied to interacting particle
systems. Widening the view to general partially ordered sets leads to a number
of new insights. In particular, we will see that:

\begin{itemize}
\item[I] Siegmund's duality for monotone processes taking values in totally
  ordered spaces and Harris' duality for additive spin systems are based on the
  same principle. Indeed, both are special cases of a general duality for
  additive processes taking values in a lattice (in the order-theoretic
  meaning of the word).
\item[II] Additive processes taking values in a general lattice (in our new,
  more general formulation) can always be constructed in terms of a
  percolation substructure. If the lattice is moreover distributive, then the
  process and its dual can be represented in the same percolation substructure
  and the duality has a graphical interpretation.
\item[III] Gray's duality for monotone spin systems can be generalized to
  monotonically representable Markov processes taking values in general
  partially ordered spaces. Here, ``monotonically representable'' is a
  somewhat stronger concept than monotonicity, as discovered by \cite{FM01}
  and independently by D.A.~Ross (unpublished).
\end{itemize}

Regarding Points~I and II, we note that our aproach not only gives a
  unified treatment of existing dualities, but also allows the systematic
  derivation of new dualities. In the context of interacting particle systems,
  all work so far has concentrated on local state spaces of the form $\{0,1\}$
  (additive systems) and $\{0,1,2\}$ (Krone's duality -see
  Section~\ref{S:Krone}), while our approach allows this to be any finite
  lattice (even local state spaces such as $\{0,1,\ldots,q\}$ with
  $q\geq 3$ have so far not been investigated). Point~II also allows to
  systematically provide a percolation picture for such processes, which is a
  convenient tool in their analysis.

The rest of the paper is organized as follows. In Section~\ref{S:main}, we
state our main results. In Subsections~\ref{S:pathdual}--\ref{S:subspace}, we
develop a general strategy for finding pathwise duals, based on invariant
subspaces, and determine the invariant subspaces associated with additive and
monotone maps, respectively, that form the basis for our pathwise dualities. In
Subsections~\ref{S:addual} and \ref{S:mondual} we then state our main results,
which are a general duality for additively representable Markov processes
(Theorem~\ref{T:addual}, corresponding to Point~I above), and general duality
for monotonically representable Markov processes (Theorem~\ref{T:mondual},
corresponding to Point~III above). In Subsection~\ref{S:percol}, we
present the results from Point~II above and in Subsection~\ref{S:open}
we mention some open problems.

Section~\ref{S:examples} is dedicated to examples. In
Subsections~\ref{S:Siegmund} and \ref{S:addpart}, we show that Siegmund's
duality for monotone processes taking values in totally ordered spaces and
Harris' duality for additive spin systems, respectively, are both special
cases of the duality for additively representable systems introduced in
Subsection~\ref{S:addual}. In Subsection~\ref{S:Krone}, we show that a nice
duality for the two-stage contact process discovered by Krone \cite{Kro99}
also fits into this general scheme. In Subsection~\ref{S:Gray}, we show that
Gray's \cite{Gra86} duality for monotone spin systems is a special case of the
general duality for monotonically representable Markov processes from
Subsection~\ref{S:mondual}. In Subsection~\ref{S:coop}, finally, we apply the
duality of Subsection~\ref{S:mondual} to derive a new duality, announced in
\cite{SS14}, for particle systems with cooperative branching and coalescence.

Section~\ref{S:proofs} is dedicated to proofs. Up to this point, all results
are for processes with finite state spaces only. In Section~\ref{S:infprod},
we generalize our results to infinite product spaces, which are the typical
setting of interacting particle systems. We conclude the paper with a short
appendix recalling some basic facts from order theory as well as a list of
notation.
%used in Section~\ref{S:pathdual}.

\section{Main results}\label{S:main}

\subsection{Duality and pathwise duality}\label{S:pathdual}

We recall that a (time-homogeneous) Markov process with measurable state
space $S$ is a stochastic process $X=(X_t)_{t\geq 0}$ taking values in $S$,
defined on an underlying probability space $(\Om,\Fi,\P)$ with expectation
denoted by $\E$, such that
\be
\label{Markov0}
\E\big[f(X_u)\,\big|\,(X_s)_{0\leq s\leq t}\big]=P_{u-t}f(X_t)\quad{\rm a.s.}
\quad(0\leq t\leq u),
\ee
for each bounded measurable $f:S\to\R$, where $(P_t)_{t\geq 0}$ is a
collection of probability kernels called the transition probabilities
of $X$, and $P_tf(x):=\int P_t(x,\di y)f(y)$.

Let $X=(X_t)_{t\geq 0}$ and $Y=(Y_t)_{t\geq 0}$ be Markov processes with state
spaces $S$ and $T$, respectively, and let $\psi:S\times T\to\R$ be a measurable
function. Then one says that $X$ and $Y$ are \emph{dual} to each other with
respect to the \emph{duality function} $\psi$, if
\be\label{dual}
\E\big[\psi(X_t,Y_0)\big]=\E\big[\psi(X_0,Y_t)\big]\qquad(t\geq 0)
\ee
for arbitrary deterministic initial states $X_0$ and $Y_0$. If (\ref{dual})
holds for deterministic initial states, then it also holds for random initial
states provided $X_t$ is independent of $Y_0$, $X_0$ is independent of $Y_t$,
and the integrals are well-defined. Possibly the first place where duality of
Markov processes is defined in this general way is \cite{CR84} (see also
\cite[Def.~II.3.1]{Lig85}), although the term duality was used much earlier
for specific duality functions. See also \cite[Section~5]{DF90} for a good
overview on the early literature on duality, as well as the relation between
duality and intertwining of Markov processes (called ``algebraic duality'' in
\cite{DF90}); see also \cite{HM11} for more on this subject.

%(\Xb_{s,t})_{s\leq t}$ is a stochastic flow. This is closely related to the
%formalism of random dynamical systems, which is a bit different, however.

It is often possible to construct Markov processes by means of a stochastic
flow. Let $S$ be a metrizable space and let $(\Xb_{s,t})_{s\leq t}$ ($s,t \in \R$) be a
collection of random maps $\Xb_{s,t}:S\to S$, such that almost surely, for
each $x\in S$, the value $\Xb_{s,t}(x)$ as a function of both $s$ and $t$ is
cadlag, i.e., right-continuous with left limits (denoted by $\Xb_{s-,t}$,
$\Xb_{s,t-}$, etc.).
\footnote{More formally, we require that there exist four functions
  $\Xb_{s,t}$, $\Xb_{s-,t}$, $\Xb_{s,t-}$, and $\Xb_{s-,t-}$ such that
  $\Xb_{s-,t}=\lim_{u\up s}\Xb_{u,t}$ and $\Xb_{s,t}=\lim_{u\down s}\Xb_{u-,t}$,
  and likewise with $t$ replaced by $t-$, and with the roles of the first and
  second variable interchanged.} We call $(\Xb_{s,t})_{s\leq t}$ a
\emph{stochastic flow} if 
\begin{enumerate}
\item $\Xb_{s,t}=\Xb_{s-,t}=\Xb_{s,t-}=\Xb_{s-,t-}$ a.s.\ for deterministic
  $s\leq t$.
\item $\Xb_{s,s}$ is the identity map and $\Xb_{t,u}\circ\Xb_{s,t}=\Xb_{s,u}$
  $(s\leq t\leq u)$.
\end{enumerate}
We say that $(\Xb_{s,t})_{s\leq t}$ has \emph{independent increments} if
\begin{enumerate}\addtocounter{enumi}{2}
\item $\Xb_{t_0,t_1},\ldots,\Xb_{t_{n-1},t_n}$ are independent for any
  $t_0<\cdots<t_n$.
\end{enumerate}
If $(\Xb_{s,t})_{s\leq t}$ is a stochastic flow with independent increments
and $X_0$ is an $S$-valued random variable, independent of $(\Xb_{s,t})_{s\leq
  t}$, then for any $s\in\R$, setting
\be\label{flowMark}
X_t:=\Xb_{s,s+t}(X_0)\qquad(t\geq 0)
\ee
defines a Markov process $(X_t)_{t\geq 0}$ with cadlag sample paths. Many
well-known Markov processes can be constructed in this way. Examples are
interacting particle systems that can be constructed from Poisson processes
that form their graphical representation, or diffusion processes that can be
constructed as strong solutions to a stochastic differential equation,
relative to a (multidimensional) Brownian motion.

Let $(\Xb_{s,t})_{s\leq t}$ and $(\Yb_{s,t})_{s\leq t}$ be stochastic flows
with independent increments, acting on metrizable spaces $S$ and $T$,
respectively, and let $\psi:S\times T\to\R$ be a
function. Then we will say that $(\Xb_{s,t})_{s\leq t}$ and
$(\Yb_{s,t})_{s\leq t}$ are \emph{dual} to each other with respect to the
\emph{duality function} $\psi$, if
\begin{enumerate}
\item
  $(\Xb_{t_0,t_1},\Yb_{-t_1,-t_0}),\ldots,(\Xb_{t_{n-1},t_n},\Yb_{-t_n,-t_{n-1}})$
  are independent for any $t_0<\cdots<t_n$.
\item For each $x\in S$, $y\in T$, and $s\leq u$, the function
  $[s,u]\ni t\mapsto\psi\big(\Xb_{s,t-}(x),\Yb_{-u,-t}(y)\big)$ is a.s.\ constant.
\end{enumerate}
Note that by our convention that stochastic flows are right-continuous, the
map $t\mapsto\Yb_{-u,-t}(y)$ is left-continuous; this is why in (ii) one needs
the left-continuous modification $\Xb_{s,t-}$ to get a sensible
definition.\footnote{The alternative is to take $\Yb_{-u,(-t)-}(y)$ which is
  right-continuous but notationally cumbersome.}
Condition~(i) implies in particular that $(\Xb_{s,t})_{s\leq t}$ and
$(\Yb_{s,t})_{s\leq t}$ have independent increments but also says that
$\Xb_{s,t}$ and $\Yb_{-t,-s}$ use, in a sense, the same randomness, i.e.,
the direction of time for the second flow is reversed with respect to the
first flow. Setting $t=s,u$ in Condition~(ii) shows that
\be
\label{dual2}
\psi\big(x,\Yb_{-u,-s}(y)\big)=\psi\big(\Xb_{s,u}(x),y\big)\quad{\rm a.s.},
\ee
which in view of (\ref{flowMark}) implies that the Markov processes $X$ and
$Y$ associated with the stochastic flows $(\Xb_{s,t})_{s\leq t}$ and
$(\Yb_{s,t})_{s\leq t}$ are dual.

If two Markov processes $X$ and $Y$ can be constructed from stochastic flows
that are dual in the sense defined above, then, loosely following terminology
introduced in \cite{JK14}, we say that $X$ and $Y$ are \emph{pathwise dual} to
each other. Although this is a much stronger concept than duality, a
remarkably large number of well-known Markov process dualities (though not
all) can be realized as pathwise dualities. In \cite{HA07,Swa06a,Swa06b}, it
is moreover demonstrated that several well-known dualities for diffusions can
be understood as continuum limits of pathwise dualities.

Sometimes (such as in \cite{AS05}) it may be useful to consider the case where
the equality in (\ref{dual}) is replaced by a $\leq$ (respectively $\geq$)
sign. If this is the case, then we say that $X$ is a \emph{subdual}
(respectively \emph{superdual}) of $Y$ (with duality function
$\psi$). Likewise, if the function $[s,u]\ni
t\mapsto\psi\big(\Xb_{s,t-}(x),\Yb_{-u,-t}(y)\big)$ is a.s.\ nonincreasing
(respectively nondecreasing), then we speak of \emph{pathwise subduality}
(respectively \emph{superduality}).

\subsection{Random mapping representations}\label{S:randmap}

Starting here, throughout the remainder of Section~\ref{S:main} as well as
Sections~\ref{S:examples} and \ref{S:proofs}, we will only be concerned with
(continuous-time) Markov processes with \emph{finite} state spaces.  Dualities
for infinite systems are often straightforward analogues of dualities for
finite systems (see Section~\ref{S:infprod}) or can be derived as limits of
the latter \cite{HA07,Kol11,Swa06a,Swa06b}. Concentrating on finite spaces
here thus allows us to separate the arguments related to duality from the
technicalities involved with infinite state spaces. As a start, in the present
subsection, we show how stochastic flows for continuous-time Markov processes
with finite state spaces can be constructed from Poisson noise and how this
leads to a standard way to construct pathwise duals.

For any two sets $S,T$, we let $\Fi(S,T)$ denote the space of all functions
$f:S\to T$. If $S,T$ are finite sets, then any linear operator
$A:\Fi(T,\R)\to\Fi(S,\R)$ is uniquely characterized by its matrix
$(A(x,y))_{x\in S,\ y\in T}$ through the formula
\be
Af(x)=\sum_{y\in T}A(x,y)f(y)\qquad\big(x\in S,\ f\in\Fi(T,\R)\big).
\ee
A linear operator $K:\Fi(T,\R)\to\Fi(S,\R)$ is a \emph{probability kernel from
  $S$ to $T$} if and only if 
\be
K(x,y)\geq 0\quand\sum_{z\in T}K(x,z)=1\qquad(x\in S,\ y\in T).
\ee
If $S$ is a finite set, then a \emph{Markov generator} (sometimes called
$Q$-matrix) of a Markov process with state space $S$ is a matrix
$(G(x,y))_{x,y\in S}$ such that
\be\label{Gdef}
G(x,y)\geq 0\quad(x\neq y)\quand\sum_{y\in S}G(x,y)=0\quad(x\in S).
\ee
It is well-known that each such Markov generator defines probability kernels
$(P_t)_{t\geq 0}$ on $S$ through the formula
\be\label{Ptdef}
P_t:=\ex{tG}=\sum_{n=0}^\infty\frac{1}{n!}t^nG^n\qquad(t\geq 0).
\ee
The $(P_t)_{t\geq 0}$ are the \emph{ transition kernels} of a Markov process
$X=(X_t)_{t\geq 0}$ with cadlag sample paths, and the Markov property
(\ref{Markov0}) can be reformulated in matrix notation as
\be\label{Markov}
\P[X_u=y\,|\,(X_s)_{0\leq s\leq t}]=P_{u-t}(X_t,y)\quad{\rm a.s.}
\qquad(0\leq t\leq u,\ y\in S).
\ee
We wish to construct the Markov process $X$ from a stochastic flow. To this
aim, we write the generator in the form
\be\label{GPois}
Gf(x)=\sum_{m\in\Gi}r_m\big(f(m(x))-f(x)\big)\qquad(x\in S, f \in \Fi(S,\R)),
\ee
where $\Gi\sub\Fi(S,S)$ is a set whose elements are maps $m:S\to S$ and
$(r_m)_{m\in\Gi}$ are nonnegative constants. Loosely following terminology
from \cite{LPW09}, we call (\ref{GPois}) a \emph{random mapping
  representation} for $G$. It is not hard to see that each Markov generator
(on a finite set $S$) can be written as in (\ref{GPois}). Random mapping
representations are far from unique, although in practical situations,
there are usually only a handful which one would deem ``natural''.

Given a random mapping representation for its generator, there is a standard
way to construct a stochastic flow for a Markov process. Let $\De$ be a
Poisson point subset of $\Gi\times\R=\{(m,t):m\in\Gi,\ t\in\R\}$ with local
intensity $r_m\di t$, where $\di t$ denotes Lebesgue measure. For $s\leq u$,
set $\De_{s,u}:=\De\cap(\Gi\times(s,u])$. Then the $\De_{s,u}$ are
a.s.\ finite sets and, because Lebesgue measure is nonatomic, two distinct
points $(m,t),(m',t')\in\De_{s,u}$ have a.s.\ different time coordinates
$t\neq t'$. Using this, we can unambiguously define random maps
$\Xb_{s,t}:S\to S$ $(s\leq t)$ by
\be\ba{l}\label{Phidef}
\dis\Xb_{s,t}:=m_n\circ\cdots\circ m_1\\[5pt]
\dis\quad\mbox{if}\quad
\De_{s,t}=\{(m_1,t_1),\ldots,(m_n,t_n)\},
\quad t_1<\cdots<t_n,
\ec
with the convention that $\Xb_{s,t}(x)=x$ if $\De_{s,t}=\emptyset$. We also
define $\De_{s,u-}:=\De\cap(\Gi\times(s,u))$, and define $\Xb_{s,t-}$
correspondingly, and similarly for $\Xb_{s-,t}$, $\Xb_{s-,t-}$. The following
lemma is well-known (see, e.g., \cite[Prop~2.5]{Swa16}).

\bl[Stochastic flow]\label{L:Markrep}
The random maps $(\Xb_{s,t})_{s\leq t}$ form a stochastic flow with
independent increments, as defined in Subsection~\ref{S:pathdual},
and the Markov process with generator $G$ can be constructed in terms of
$(\Xb_{s,t})_{s\leq t}$ as in (\ref{flowMark}).
\el

Let $S$ and $T$ be sets and let $\psi:S\times T\to\R$ be a function. Then we
say that two maps $m:S\to S$ and $\hat m:T\to T$ are \emph{dual} with respect
to the \emph{duality function} $\psi$ if
\be\label{mdual}
\psi\big(m(x),y\big)=\psi\big(x,\hat m(y)\big)
\qquad(x\in S,\ y\in T).
\ee
If in (\ref{mdual}) the equality is replaced by $\leq$ (resp.\ $\geq$) then we
also say that $\hat m$ is \emph{subdual} (resp.\ \emph{superdual}) to $m$ with
respect to the duality function $\psi$. We remark that for a given $\psi$ and
$m$ there may be many dual maps $\hat m$. However, the dual map $\hat m$ is
unique if the dual function has the property that $\psi(x,y)=\psi(x,z)$ for
all $x\in S$ implies $y=z$, for general $y,z\in T$.

Now imagine that $S$ and $T$ are
finite sets, that $G$ is the generator of a Markov process in $S$, and that
for a given random mapping representation as in (\ref{GPois}), all maps
$m\in\Gi$ have a dual $\hat m$ with respect to $\psi$.
 %(In general, such a dual map need not be unique, but we choose one and denote it by $\hat m$.)
Then we claim that the Markov process $Y$ with state space $T$ and generator
\be
\label{HPois}
Hf(y):=\sum_{m\in\Gi}r_m\big(f(\hat m(y))-f(y)\big)
\qquad(y\in T,\ f\in\Fi(T,\R))
\ee
is pathwise dual to $X$. To see this, we define
\be\label{hatDe}
\hat\De:=\{(\hat m,-t):(m,t)\in\De\},
\ee
which is a Poisson point set on $\hat\Gi\times\R$ with $\hat\Gi:=\{\hat
m:m\in\Gi\}$. We use this Poisson point set to define random maps
$(\Yb_{s,t})_{s\leq t}$ in the same way as in (\ref{Phidef}). The proof of the
following proposition, which can just as well be formulated in terms of sub-
or superduality, is entirely straightforward; for completeness, we give the
main steps in Section~\ref{S:Markdual}.

\bp[Pathwise duality]\label{P:pathdual}
Let $X$ and $Y$ be Markov processes with generators $G$ and $H$ of the
  form (\ref{GPois}) and (\ref{HPois}), respectively, where for each
  $m\in\Gi$, the map $\hat m$ is a dual of $m$ in the sense of (\ref{mdual}).
  Construct stochastic flows $(\Xb_{s,t})_{s\leq t}$ and $(\Yb_{s,t})_{s\leq
    t}$ for these Markov processes as above. Then, almost surely, the map
$\Xb_{s-,t-}$ is dual to $\Yb_{-t,-s}$ for each $s\leq t$. Moreover, the
stochastic flow $(\Yb_{s,t})_{s\leq t}$ is dual to $(\Xb_{s,t})_{s\leq t}$
with respect to the duality function $\psi$, in the sense defined in
Subsection~\ref{S:pathdual}, and the Markov processes $X$ and $Y$ are pathwise
dual w.r.t.\ $\psi$.
\ep

By grace of Proposition~\ref{P:pathdual}, in order to prove that a given
Markov process $X$ with generator $G$ has a pathwise dual $Y$ with respect to
a certain duality function $\psi$, and in order to explicitly construct such a
dual, it suffices to show that $G$ has a random mapping representation
(\ref{GPois}) such that each map $m\in\Gi$ has a dual w.r.t.\ $\psi$. In view
of this, much of our paper will be devoted to showing that certain maps have
duals with respect to particular duality functions. Once we have shown this
for a suitable class of maps, it is clear how to construct pathwise duals for
Markov processes whose generators are representable in such maps.

In relation to this, we adopt the following definition. Let $S$ be a finite
set and let $\Gi\sub\Fi(S,S)$ be a set whose elements are maps $m:S\to S$.  We
say that a Markov generator $G$ is \emph{representable in $\Gi$} if $G$ can be
written in the form (\ref{GPois}) for nonnegative constants $(r_m)_{m\in\Gi}$.

Similar definitions apply to probability kernels. Let $S$ and $T$ be finite
sets. We let $\Ki(S,T)$ denote the space of all probability kernels from $S$
to $T$. If $M$ is a random variable taking values in the space $\Fi(S,T)$, then
\be\label{randmap}
K(x,y):=\P[M(x)=y]\qquad(x\in S,\ y\in T)
\ee
defines a probability kernel $K\in\Ki(S,T)$. If a probability kernel $K$ is
written in the form (\ref{randmap}), then, borrowing terminology from
\cite{LPW09}, we call (\ref{randmap}) a \emph{random mapping representation}
for $K$. For any set $\Gi\sub\Fi(S,T)$ of maps $m:S\to T$, we say that a
probability kernel $K\in\Ki(S,T)$ is \emph{representable in $\Gi$} if there
exists a $\Gi$-valued random variable $M$ such that (\ref{randmap}) holds. It
is easy to see that each probability kernel from $S$ to $T$ is representable
in $\Fi(S,T)$ but the representation is again not unique.

We observe that if $(\Xb_{s,t})_{s\leq t}$ is a stochastic flow with
independent increments, then the transition kernels of the associated Markov
process (in the sense of (\ref{flowMark})) are given by
\be\label{PXb}
P_t(x,y)=\P[\Xb_{0,t}(x)=y]\qquad(t\geq 0,\ x,y\in S),
\ee
and this formula gives a random mapping representation for $P_t$. In view of
this and (\ref{Phidef}), the following lemma, which will be proved in Section
\ref{S:Markdual}, should not come as a surprise.

\bl[Representability of generators]\label{L:genrep}
Let $S$ be a finite set, let $\Gi$ be a set whose elements are maps $m:S\to
S$, and let $X$ be a Markov process in $S$ with generator $G$ and transition
kernels $(P_t)_{t\geq 0}$. Assume that $\Gi$ is closed under composition
and contains the identity map. Then
the following statements are equivalent:
\begin{enumerate}
\item $G$ can be represented in $\Gi$.
\item $P_t$ can be represented in $\Gi$ for all $t\geq 0$.
\end{enumerate}
\el

\subsection{Invariant subspaces}\label{S:subspace}

As explained previously we will due to Proposition~\ref{P:pathdual} be
interested in Markov processes whose generators can be represented in maps
that have a dual map with respect to a suitable duality function. In the
present subsection, we describe a general strategy for finding such maps, and
for choosing the duality function. The following simple observation shows that
in fact, each map is dual to its inverse image map. For any set $S$, we let
$\Pc(S):=\{A:A\sub S\}$ denote the set of all subsets of $S$.  Below, and in
what follows, we use the notation $1_{\{\ldots\}}:=1$ if the statement
$\ldots$ holds, and $1_{\{\ldots\}}:=0$ otherwise. So, for example,
$x\mapsto1_{\{x\leq 1\}}$ is the indicator function of the set $\{x:x\leq 1\}$
and $(x,y)\mapsto1_{\{x\leq y\}}$ is the indicator function of the set
$\{(x,y):x\leq y\}$.

\bl[Inverse image map]\label{L:invdual}
Let $S$ be a set, let $m:S\to S$ be a map, and let $m^{-1}:\Pc(S)\to\Pc(S)$ be
the inverse image map defined as $m^{-1}(A):=\{x\in S:m(x)\in A\}$. Then
$m^{-1}$ is dual to $m$ with respect to the duality function
$\psi(x,A):=1_{\{x\in A\}}$.
\el

Combining this with Proposition~\ref{P:pathdual}, we see
  that if the generator $G$ of a Markov process $X$
with state space $S$ has a random mapping representation of the form
(\ref{GPois}), then the Markov process $Y$ with state space $\Pc(S)$ and
generator $H$ given by
\be\label{Hdef}
Hf(A)=\sum_{m\in\Gi}r_m\big(f(m^{-1}(A))-f(A)\big)\qquad\big(A\in\Pc(S)\big)
\ee
is a pathwise dual of $X$ with respect to the duality function $\psi$ from
Lemma~\ref{L:invdual}. In practice, this dual process is not very useful since
the state space $\Pc(S)$ is very large compared to $S$. The situation is
better, however, if $\Pc(S)$ contains a smaller subspace that is
 invariant\footnote{The connection between pathwise
    duality and invariant subspaces exposed here has an analogue for
    (non-pathwise) duality. In \cite[Prop~2.2]{Moe13}, it is shown that if
    $S,T$ are finite sets and $\psi:S\times T\to\R$ is such that the
    functions $\{\psi(\,\cdot\,,y):y\in T\}$ are linearly independent, then a
    Markov process $X$ with state space $S$ has a dual w.r.t.\ to the duality
    function $\psi$ if and only if the semigroup of $X$ maps the convex hull of
    $\{\psi(\,\cdot\,,y):y\in T\}$ into itself.} under
all maps $m\in\Gi$. In the present paper, we will be interested in the
subspace of all decreasing subsets of $S$ (where $S$ is a partially ordered
set) and the subspace of all principal ideals of $S$ (where $S$ is a lattice),
which leads to a duality for Markov processes that are representable in
monotone maps and additive maps, respectively.

So let $S$ now be a partially ordered set and let us briefly introduce some 
definitions and recall some basic facts for partially ordered sets.
For any set $A\sub S$ we define
$A^\up:=\{x\in S:x\geq y\mbox{ for some }y\in A\}\supset A$. We say that $A$ is
\emph{increasing} if $A^\up\sub A$. We define $A^\down$ and \emph{decreasing}
sets in the same way for the reversed order. A nonempty increasing set $A$
such that for every $x,y\in A$ there exists a $z\in A$ with $z\leq x,y$ is
called a \emph{filter} and a nonempty decreasing set $A$ such that for every
$x,y\in A$ there exists a $z\in A$ with $x,y\leq z$ is called an
\emph{ideal}. A \emph{principal filter} is a filter that contains a minimal
element and a \emph{principal ideal} is an ideal that contains a maximal
element. Equivalently, principal filters are sets of the form $A=\{z\}^\up$
and principal ideals are sets of the form $A=\{z\}^\down$, for some $z\in S$.
A finite filter or ideal is always principal. We will use the notation
\bc\label{!dec}
\dis\Pc_{\rm inc}(S)&:=&\dis\{A\sub S:A\mbox{ is increasing}\},\\[5pt]
\dis\Pc_{\rm !inc}(S)&:=&\dis\{A\sub S:A\mbox{ is a principal filter}\},\\[5pt]
\dis\Pc_{\rm dec}(S)&:=&\dis\{A\sub S:A\mbox{ is decreasing}\},\\[5pt]
\dis\Pc_{\rm !dec}(S)&:=&\dis\{A\sub S:A\mbox{ is a principal ideal}\}.
\ec

A \emph{lattice} is a partially ordered set for which the sets of principal
filters and principal ideals are closed under finite intersections.
Equivalently, this says that for every $x,y\in S$ there exist (necessarily
unique) elements $x\vee y\in S$ and $x\wedge y\in S$ called the
\emph{supremum} or \emph{join} and \emph{infimum} or \emph{meet} of $x$ and
$y$, respectively, such that
\be\label{supdef}
\{x\}^\up\cap\{y\}^\up=:\{x\vee y\}^\up
\quand
\{x\}^\down\cap\{y\}^\down=:\{x\wedge y\}^\down.
\ee
A \emph{join-semilattice} (respectively \emph{meet-semilattice}) is a partially
ordered set in which $x\vee y$ (respectively $x\wedge y$) are well-defined. A
partially ordered set $S$ is \emph{bounded from below} if it contains an
(obviously unique) element, usually denoted by 0, such that $0\leq x$ for all
$x\in S$. Boundedness from above is defined analogously and the (obviously
unique) upper bound is often denoted by 1. Finite lattices are always bounded
from below and above.

Let $S$ and $T$ be partially ordered sets. By definition, a map $m:S\to T$
is \emph{monotone} if it is a $\leq$-homomorphism, i.e.,
\be
\label{mondef}
x\leq y\quad\mbox{implies}\quad m(x)\leq m(y)\quad(x,y\in S).
\ee
We denote the set of all monotone maps $m:S\to T$ by $\Fi_{\rm mon}(S,T).$
If $S$ and $T$ are join-semilattices that are bounded from below, then,
generalizing terminology from the theory of interacting particle systems, we
will call a map $m:S\to T$ \emph{additive} if it is a
$(0,\vee)$-homomorphism, i.e.,
\be\label{addef}
m(0)=0\quand m(x\vee y)=m(x)\vee m(y)\qquad(x,y\in S).
\ee
From these definitions it is easy to see that any additive map is also
monotone.

Our main duality results concern Markov processes that are \emph{monotonically
  representable} by which we mean that they possess a random mapping
representation as in (\ref{GPois}) with ${\cal G}= \Fi_{\rm mon}.$
Respectively, we are considering \emph{additively representable} processes for
which ${\cal G}$ contains only additive functions.

The basis of all these dualities is the following simple lemma, which says
that a map $m$ is monotone, respectively additive, if and only if its inverse
image map leaves the subspaces of all decreasing subsets, respectively ideals,
invariant. Note that since in a finite lattice, all ideals are principal,
part~(ii) says in particular that if $S$ and $T$ are finite join-semilattices
that are bounded from below, then $m:S\to T$ is additive if and only if
$m^{-1}(A)\in\Pc_{\rm !dec}(S)$ for all $A\in\Pc_{\rm !dec}(T)$. The proof of
Lemma~\ref{L:monadd} will be given in Section~\ref{S:monad}.

\bl[Monotone and additive maps]
{\rm (i)} Let $S$ and $T$ be partially ordered sets and let $m:S\to T$ be a
map. Then $m$ is monotone if and only if\label{L:monadd}
\begin{quote}
$m^{-1}(A)\in\Pc_{\rm dec}(S)$ for all $A\in\Pc_{\rm dec}(T)$.
\end{quote}
{\rm (ii)} If $S$ and $T$ are join-semilattices that are bounded from below,
then $m$ is additive if and only if
\begin{quote}
$m^{-1}(A) \sub S$ is an ideal whenever $A\sub T$ is an ideal.
\end{quote}
\el

\subsection{Additive systems duality}\label{S:addual}

In the present section, we show how the fact that the inverse image of an
additive map leaves the subspace of principal ideals invariant, leads in a
natural way to a duality for Markov processes that are additively
representable.

We start with an abstract definition of a dual
for any partially ordered set $S.$
Namely, a \emph{dual} of $S$ is a
partially ordered set $S'$ together with a bijection $S\ni x\mapsto x'\in S'$
such that
\be\label{duspa}
x\leq y\quad\mbox{if and only if}\quad x'\geq y'.
\ee
Two canonical ways to construct such a dual are as follows. \emph{Example~1:}
For any partially ordered set $S$, we may take $S':=S$ but equipped with the
reversed order, and $x\mapsto x'$ the identity map. \emph{Example~2:} If $\La$
is a set and $S\sub\Pc(\La)$ is a set of subsets of $\La$, equipped with the
partial order of inclusion, then we may take for $x':=\La\beh x$ the
complement of $x$ and $S':=\{x':x\in S\}$.

Returning to the abstract definition, it is easy to see that all duals of a
partially ordered set are naturally isomorphic and that the original partially
ordered set is in a natural way the dual of its dual, which motivates us to
write $x''=x$. If $S$ is bounded from below, then $S'$ is bounded from above
and $0'=1$. We define a function $S\times S'\ni(x,y)\mapsto\li
x,y\re\in\{0,1\}$ by
\be\label{tetSS}
\li x,y\re:=1_{\txt\{x\leq y'\}}=1_{\txt\{y\leq x'\}}\qquad(x\in S,\ y\in S').
\ee
Note that $x,y'\in S$ while $y,x'\in S'$, so the two signs $\leq$ in this
formula refer to the partial orders on $S$ and $S'$, respectively. Since $S$
is the dual of $S'$, the formal definition of $\li\,\cdot\,,\,\cdot\,\re$ is
symmetric in the sense that $\li y,x\re=\li x,y\re$ for all $y\in S'$ and
$x\in S''=S$. In the concrete examples of a dual space given above,
$\li\,\cdot\,,\,\cdot\,\re$ has the following meanings. In Example~1, $\li
x,y\re=1_{\{x\leq y\}}$, while in Example~2, $\li x,y\re=1_{\{x\cap
  y=\emptyset\}}$.

Since our aim is to show how additive systems duality arises naturally from
Lemma~\ref{L:monadd}, we will now let $S$ be a lattice and prove the following lemma here on the spot.

\bl[Duals of additive maps]\label{L:addual}
Let $S$ be a finite lattice and let $S'$ be a dual of $S$ in the sense defined
above. Then a map $m:S\to S$ is additive if and only if there exists a
(necessarily unique) map $m':S'\to S'$ that is dual to $m$ with respect to the
duality function $\psi(x,y):=\li x,y\re$ from (\ref{tetSS}). This dual map
$m'$, if it exists, is also additive.
\el
\begin{Proof}
By Lemma~\ref{L:monadd}, $m$ is additive if and only if $m^{-1}(A)\in\Pc_{\rm
  !dec}(S)$ for all $A\in\Pc_{\rm !dec}(S)$. Since each element $A\in\Pc_{\rm
  !dec}(S)$ can be written in the form $A=\{x\in S:x\leq y'\}$ for some unique
$y\in S'$, $m$ is additive if and only if for each $y\in S'$ there exists a
(necessarily unique) element $m'(y)\in S'$ such that
\be
m^{-1}\big(\{z\in S:z\leq y'\}\big)=\{x\in S:x\leq(m'(y))'\},
\ee
i.e., there exists a (necessarily unique) map $m':S'\to S'$ such that
\be
m(x)\leq y'\quad\mbox{if and only if}\quad x\leq(m'(y))'
\qquad(x\in S,\ y\in S'),
\ee 
which says that $\li m(x),y\re=\li x,m'(y)\re$, i.e., $m'$ is dual to $m$ with
respect to the duality function from (\ref{tetSS}).

Since conversely $m$ is dual to $m'$ with respect to the duality function from
(\ref{tetSS}), by what we have just proved, $m'$ (having a dual with respect
to this duality function) must be additive.
\end{Proof}

Combining Proposition~\ref{P:pathdual} and Lemma~\ref{L:addual}, we can write
down our first nontrivial result.

\bt[Additive systems duality]\label{T:addual}
Let $S$ be a finite lattice and let $X$ be a Markov process in $S$ whose
generator has a random mapping representation of the form
\be\label{GPoisadd}
Gf(x)=\sum_{m\in\Gi}r_m\big(f(m(x))-f(x)\big)\qquad(x\in S),
\ee
where all maps $m\in\Gi$ are additive. Then the Markov process $Y$ in $S'$
with generator
\be\label{GPoisY}
Hf(y):=\sum_{m\in\Gi}r_m\big(f(m'(y))-f(y)\big)\qquad(y\in S')
\ee
is pathwise dual to $X$ with respect to the duality function $\psi(x,y):=\li
x,y\re$ from (\ref{tetSS}).
\et

In Subsections~\ref{S:Siegmund} and \ref{S:addpart} below, we will see that
Theorem~\ref{T:addual} contains both Siegmund's duality and the well-known
duality of additive interacting particle systems as special cases.  Moreover,
in Section~\ref{S:Krone}, we will see that a duality for the two-stage contact
process discovered by Krone \cite{Kro99} also fits into this general scheme.
We note that the similarity between Siegmund's duality and additive systems
duality was already observed in \cite{CS85}, who argued that both dualities
could be formulated in terms of processes taking values in sets of subsets. In
Section~\ref{S:percol} we will discuss in detail when additive processes and
their duals can be represented as set-valued processes.

Finally, we observe that duality (not necessarily pathwise) with respect to
the the duality function $\psi(x,y):=\li x,y\re$ immediately implies that the
underlying processes are monotone Markov processes in the traditional
sense. This appears to be a well-known fact even in relatively general
settings and follows since due to the duality for any $y \in S'$ the function
$x \mapsto \P(X_t > y'|X_0=x)$ is a monotone function. In other words, for the
monotone functions $f_{y'}(x) = 1_{\{x>y'\}}$ we have that $P_t f_{y'}$ is
again monotone where $(P_t)_{t \geq 0}$ are the transition kernels of $X.$
This implies in turn the monotonicity, see (\ref{monK}) below, as any monotone
function can be written as a linear combination of the functions $f_{y'}, y
\in S'$.

The question remains whether additive representability is itself a necessary
condition for the existence of a (again not necessarily pathwise) dual with
respect to the duality function $\psi(x,y):=\li x,y\re.$ Some sufficient
conditions to this respect for a Markov process with values in $\R^d$
(equipped with the product order) are given in \cite[Thms~4, 8, 9]{Lee13}.

\subsection{Monotone systems duality}\label{S:mondual}

In the present section, we show how the fact that the inverse image of a
monotone map leaves the subspace of decreasing subsets invariant can be used
to construct pathwise duals of general monotonically representable Markov
processes.  These dual processes are more complicated than in the case of
additively representable processes, but we will stay as close as possible to
the formalism of the previous subsection, so that Theorem~\ref{T:addual} will
be a special case of a more general theorem to be formulated here.

Theorem~\ref{T:addual} is based on the fact that if $m$ is an additive map,
then its inverse image maps the space of principal ideals into itself, and
each principal ideal $A\in\Pc_{\rm !dec}(S)$ can be encoded in terms of the
unique element $y\in S'$ such that $A=\{y'\}^\down$. In our present setting,
we will encode a general decreasing set $A\in\Pc_{\rm dec}(S)$ in terms of a
set $B\sub S'$ such that $A=\{y':y\in B\}^\down$. This means that we will use
the duality function
\be\label{phidef}
\phi(x,B):=1_{\txt\{x\leq y'\mbox{ for some }y\in B\}}
\qquad(x\in S,\ B\sub S').
\ee
For a given $A\in\Pc_{\rm dec}(S)$, there are usually more ways to choose a
set $B\sub S'$ such that $A=\{y':y\in B\}^\down$ and as a result we see that
for a given monotone map $m$ there are at least two natural ways to
define a dual map with respect to the duality function $\phi$ from
(\ref{phidef}).

As before, we assume that $S$ is a finite partially ordered set and we let $S'$
(together with the map $x\mapsto x'$) denote a dual of $S$ in the sense
defined in (\ref{duspa}). Contrary to the latter part of the previous subsection, we no longer
assume that $S$ is a lattice. For any set $A\sub S$ we write $A':=\{x':x\in
A\}$ and we let
\be\label{Amax}
A_{\rm max}:=\{x\in A:x\mbox{ is a maximal element of }A\}
=\{x\in A:\not\!\exists   y \in A, y \neq x\mbox{ s.t.\ }x\leq y\}
\ee
denote the set of maximal elements of $A$. Similarly, we let $A_{\rm min}$
denote the set of minimal elements of a set $A$. For any monotone $m:S\to S$,
we can uniquely define maps $m^\dgg:\Pc(S')\to\Pc(S')$ and
$m^\ast:\Pc(S')\to\Pc(S')$ by
\be\label{mast}
m^\dgg(B)':=(m^{-1}({B'}^\down))_{\rm max}
\quand
m^\ast(B)':=\bigcup_{x\in B}(m^{-1}(\{x'\}^\down))_{\rm max}
\ee
$(B\in\Pc(S'))$. The next lemma shows that both $m^\dgg$ and $m^\ast$ are dual
to $m$ with respect to the duality function $\phi$ from (\ref{phidef}). In
addition, both $m^\dgg$ and $m^\ast$ each have a special property justifying
their definitions.

\bl[Duals of monotone maps]\label{L:dumaps}
Let $S$ be a finite partially ordered set and let $m\in\Fi_{\rm mon}(S,S)$.
Then $m^\dgg$ and $m^\ast$ are dual to $m$ with
respect to the duality function $\phi$ from (\ref{phidef}). Moreover,
\be\label{dagprop}
m^\dgg(B)=m^\dgg(B)_{\rm min}=m^\ast(B)_{\rm min}\qquad\big(B\in\Pc(S')\big),
\ee
and
\be\label{cuprop}
m^\ast(B\cup C)=m^\ast(B)\cup m^\ast(C)\qquad\big(B,C\in\Pc(S')\big).
\ee
\el

We note that if we restrict the state space for the dual process to the set of
all $B\in\Pc(S')$ such that $B=B_{\rm min}$, then the duality function $\phi$
has the property that $\phi(x,B)=\phi(x,C)$ for all $x\in S$ implies $B=C$; as
a result, the dual map is now unique and in fact given by $m^\dgg$.
Nevertheless, it is sometimes useful to choose the larger state space
$\Pc(S')$ for the dual process and in this case uniqueness of the dual map is
lost. In fact, the alternative dual map $m^\ast$ now has the nice property
(\ref{cuprop}) which sometimes makes it preferable over $m^\dgg$.

The proof of Lemma~\ref{L:dumaps} is straightforward but a bit tedious, and
for this reason we postpone it till Section~\ref{S:monproof}. Combining
Proposition~\ref{P:pathdual} and Lemma~\ref{L:dumaps}, we can write down our
second nontrivial result.\footnote{Theorem~\ref{T:mondual} shows that a
  process $X$ has a pathwise dual w.r.t.\ to the duality function $\phi$ from
  (\ref{phidef}) if $X$ is monotonically representable. This condition is also
  necessary: if a map $m$ has a dual w.r.t.\ to $\phi$, then $m^{-1}$ maps
  decreasing sets into decreasing sets and hence, by Lemma~\ref{L:monadd}, $m$
  is monotone.}

\bt[Monotone systems duality]\label{T:mondual}
Let $S$ be a finite partially ordered set and let $X$ be a Markov process in
$S$ whose generator has a random mapping representation of the form
\be\label{GPoismon}
Gf(x)=\sum_{m\in\Gi}r_m\big(f(m(x))-f(x)\big)\qquad(x\in S),
\ee
where all maps $m\in\Gi$ are monotone. Then the $\Pc(S')$-valued Markov
processes $Y^\dgg$ and $Y^\ast$ with generators
\be\ba{r@{\,}c@{\,}l}\label{HPoisY}
\dis H_\dgg f(B)&:=&\dis\sum_{m\in\Gi}r_m\big(f(m^\dgg(B))-f(B)\big),\\[5pt]
\dis H_\ast f(B)&:=&\dis\sum_{m\in\Gi}r_m\big(f(m^\ast(B))-f(B)\big),\\[5pt]
\ea
\qquad\big(B\in\Pc(S')\big)
\ee
are pathwise dual to $X$ with respect to the duality function $\phi$ from
(\ref{phidef}).
\et

We note that if a map $m:S\to S$ is monotone, then it is also monotone with
respect to the reversed order on $S$. As a result, the inverse image map
$m^{-1}$ also maps increasing subsets into increasing subsets. This naturally
leads to the duality function (compare (\ref{phidef}))
\be\label{altphi}
\ti\phi(x,B):=1_{\txt\{x\geq y'\mbox{ for some }y\in B\}}
\qquad(x\in S,\ B\sub S').
\ee
In analogy with (\ref{mast}), we may define maps $m^\circ:\Pc(S')\to\Pc(S')$ and
$m^\dott:\Pc(S')\to\Pc(S')$ 
by
\be\label{mstar}
m^\circ(B)':=(m^{-1}({B'}^\up))_{\rm min}
\quand
m^\dott(B)':=\bigcup_{x\in B}(m^{-1}(\{x'\}^\up))_{\rm min}.
\ee
As a direct consequence of Lemma~\ref{L:dumaps}, applied to the reversed order
on $S$, we obtain the following lemma.

\bl[Alternative duals of monotone maps]\label{L:altmaps}
Let $S$ be a finite partially ordered set and let $m\in\Fi_{\rm mon}(S,S)$.
Then $m^\circ$ and $m^\dott$ are dual to $m$ with
respect to the duality function $\ti\phi$ from (\ref{altphi}). Moreover,
\be\label{altprop}
m^\circ(B)=m^\circ(B)_{\rm max}=m^\dott(B)_{\rm max}\qquad\big(B\in\Pc(S')\big),
\ee
and
\be\label{cualt}
m^\dott(B\cup C)=m^\dott(B)\cup m^\dott(C)\qquad\big(B,C\in\Pc(S')\big).
\ee
\el

We will mostly focus on the maps $m^\dgg$ and $m^\ast$ from (\ref{mast}) since
these are most closely related to the additive systems duality from the
previous subsection. The next lemma, which will be proved in
Section~\ref{S:monproof}, shows that the dual processes of
Theorem~\ref{T:mondual} reduce to the dual of Theorem~\ref{T:addual} if all
maps occurring in (\ref{GPoismon}) are additive and the dual process
$Y=Y^\dgg$ or $=Y^\ast$ is started in a singleton, i.e., a state of the form
$Y_0=\{y\}$ for some $y\in S'$.

\bl[Relation between additive and monotone duals]\label{L:astac}
Let $S$ be a finite lattice and let $m:S\to S$ be additive. Let $m'$ be the
dual map from Lemma~\ref{L:addual} and let $m^\ast$ and $m^\dgg$ be as in
(\ref{mast}). Then
\be\label{astac}
m^\ast(B)=\{m'(y):y\in B\}
\quand
m^\dgg(B)=m^\ast(B)_{\rm min}
\qquad\big(B\in\Pc(S')\big).
\ee
\el

In Subsection~\ref{S:Gray}, we will show that if $X$ is a monotone spin
system, then the dual process $Y^\dott$ coincides
with Gray's \cite{Gra86} dual process. Here $Y^\dott$ denotes the process
with generator $H_\dott$ which is defined analogously to $H_\dgg,H_\ast$ in
(\ref{HPoisY}), but with the dual map $m^\dott$ from (\ref{mstar}) instead of
$m^\dgg,m^\ast$. As a more concrete application of this sort of dualities, in
Subsection~\ref{S:coop}, we use Theorem~\ref{T:mondual} to derive a new
duality, announced in \cite{SS14}, for particle systems with cooperative
branching.

Let us finish this section by relating the concepts of monotonically
representable and monotone Markov processes.  A Markov process is
traditionally called \emph{monotone} if its transition kernels $(P_t)_{t\geq
  0}$ are monotone, where, in general, for $S$ and $T$ finite partially
ordered sets, a probability kernel $K$ from $S$ to $T$ is called
\emph{monotone} if
\be\label{monK}
f\in\Fi_{\rm mon}(S,\R)\quad\mbox{implies}\quad Kf\in\Fi_{\rm mon}(S,\R),
\ee
with as before $Kf(x):=\sum_yK(x,y)f(y)$.

It turns out that being monotonically representable is a stronger concept than
being monotone.  Indeed, by Lemma \ref{L:genrep}, if the generator of a Markov
process is monotonically representable, then the same is true for its
transition kernels $(P_t)_{t\geq 0}$ and it is easy to see that each
monotonically representable probability kernel is also monotone. However, it
is known that there exist kernels that are monotone yet not monotonically
representable. See \cite[Example~1.1]{FM01} for an example where
$S=T=\Pc(\La)$ with $\La$ a set containing just two elements. (Note that what
we call monotonically representable is called \emph{realizably monotone} in
\cite{FM01}.) On the positive side, we cite the following result from
\cite{KKO77,FM01}.

\bp[Sufficient conditions for monotone representability]\label{P:monsuf}
Let $S,T$ be finite partially ordered sets and assume that at least one of the
following conditions is satisfied:
\begin{enumerate}
\item $S$ is totally ordered.
\item $T$ is totally ordered.
\end{enumerate}
Then any monotone probability kernel from $S$ to $T$ is monotonically
representable.
\ep
\begin{Proof}
The sufficiency of (i) was proved in \cite{KKO77} and also follows as a
special case of \cite[Thm~4.1]{FM01}, where it is shown that the statement
holds more generally if $S$ has a tree-like structure.
The sufficiency of (ii) is proved in \cite[Example~1.2]{FM01}.
\end{Proof}

The question how to determine whether a given probability kernel (or Markov
generator) can be represented in the set of additive maps was already
mentioned in \cite{Har78}, but we do not know about any results in this
direction. Cursory contemplation suggests that this problem is at least as
difficult as monotone representability.

\subsection{Percolation substructures}\label{S:percol}

%It seems the term Percolation substructure comes from \cite{Gri79}
%The first graphical representations for GENERAL additive systems
%occur in \cite{Har78}.

If $\La$ is a finite set, then the set $S:=\Pc(\La)$ of all subsets of $\La$,
equipped with the order of set inclusion, is a finite lattice. For lattices of
this form, there exists a simple way to characterize additive maps $m:S\to S$
and this naturally leads to a representation of additive Markov processes in
$S$, and their duals, in terms of a form of oriented percolation on
$\La\times\R$, where the second coordinate plays the role of time
\cite{Har78,Gri79}. In the present subsection, we investigate whether this
picture remains intact when $S$ is replaced by a more general lattice. For
distributive lattices, the answer is largely positive. In
Subsection~\ref{S:Krone} below, we apply the results of the present subsection
to give a percolation representation for the two-stage contact process and its
dual, which are interacting particle systems with a state space of the form
$\{0,1,2\}^\La$ \cite{Kro99}.

By definition, a lattice $S$ is \emph{distributive} if
\be\label{distrib}
x\wedge(y\vee z)=(x\wedge y)\vee(x\wedge z)\qquad(x,y,z\in S).
\ee
If $\La$ is a finite set and $S\sub\Pc(\La)$ is closed under intersections and
unions and contains $\emptyset$ and $\La$ as elements, then $S$ is a
distributive lattice with $\emptyset$ and $\La$ as lower and upper bounds. In
particular, if $\La$ is a partially ordered set, then this applies to
$S:=\Pc_{\rm dec}(\La)$. Birkhoff's representation theorem (\cite{Bir37}; see,
e.g., \cite[Thm~5.12]{DP02} for a modern reference) says that conversely, each
distributive lattice is of this form. Note that in particular, if $\La$ is
equipped with the \emph{trivial order} $i\not\leq j$ for all $i\neq j$, then
$\Pc_{\rm dec}(\La)=\Pc(\La)$.

In view of this, let $\La$ be a finite partially ordered set, and let
$S:=\Pc_{\rm dec}(\La)$ be the lattice of decreasing subsets of $\La$. Let
$\La'$ denote the set $\La$ equipped with the reversed order. Then
$S':=\Pc_{\rm dec}(\La')=\Pc_{\rm inc}(\La)$, together with the complement map
$x\mapsto x':=x^{\rm c}$, is a dual of the lattice $S$ in the sense defined in
(\ref{duspa}). The following lemma will be proved in Section~\ref{S:percproof}.

\bl[Characterization of additive maps]\label{L:madd}
There is a one-to-one correspondence between, on the one hand, additive maps
$m:S\to S$ and, on the other hand, sets $M\sub\La\times\La$ such that for all
$i,j,\ti\ii, \ti\jj \in\La$
\be\ba{rl}\label{Mprop}
{\rm(i)}&\mbox{$(i,j)\in M$ and $i\leq\ti\ii$ implies $(\ti\ii,j)\in M$,}\\[5pt]
{\rm(ii)}&\mbox{$(i,j)\in M$ and $j\geq\ti\jj$ implies $(i,\ti\jj)\in M$.}
\ec
This one-to-one correspondence comes about by identifying $M$ with the
additive map $m$ defined by
\be\label{minM}
m(x):=\{j\in\La:(i,j)\in M\mbox{ for some }i\in x\}
\qquad(x\in S).
\ee
If $m':S'\to S'$ is the dual of $m$ (in the sense of Lemma~\ref{L:addual}) and
$M'$ is the corresponding subset of $\La'\times\La'$, then $M$ and $M'$ are
related by
\be\label{Madj}
M'=\big\{(j,i):(i,j)\in M\big\}.
\ee
\el

Using Lemma~\ref{L:madd}, we can represent the stochastic flow of an additive
Markov process with values in a general distributive lattice in terms of open
paths in a ``percolation substructure'' in the sense of \cite{Gri79}. Let $X$
be a Markov process whose generator $G$ has a random mapping representation of
the form (\ref{GPois}) where all maps $m\in\Gi$ are additive maps $m:S\to S$,
with $S=\Pc_{\rm dec}(\La)$ as before. As in Section~\ref{S:randmap}, we
construct a Poisson point set $\De$ on $\Gi\times\R$ and use this to define a
stochastic flow $(\Xb_{s,t})_{s\leq t}$ as in (\ref{Phidef}). We also set
$\De':=\{(m',-t):(m,t)\in\De\}$ (compare (\ref{hatDe})) and use this Poisson
set to define the dual stochastic flow $(\Yb_{s,t})_{s\leq t}$ as in
Proposition~\ref{P:pathdual}.

Plotting space-time $\La\times\R$ with time upwards, for
each $(m,t)\in\De$, let $M$ be the set corresponding to $m$ in the sense of
Lemma~\ref{L:madd}, draw an arrow from $(i,t)$ to $(j,t)$ for each $i\neq j$
such that $(i,j)\in M$, and place a ``blocking symbol'' \block at $(i,t)$
whenever $(i,i)\not\in M$. Examples of such graphical representations are
given in Figures~\ref{fig:voter} and \ref{fig:Krone} below.
By definition, an \emph{open path} in such a graphical representation is
a cadlag function $\ga:[s,u]\to\La$ such that:
\begin{enumerate}
\item If $\ga_{t-}\neq\ga_t$ for some $t\in(s,u]$, then there is an arrow from
  $(\ga_{t-},t)$ to $(\ga_t,t)$.
\item If there is a blocking symbol at $(\ga_t,t)$ for some $t\in(s,u]$, then
  $\ga_{t-}\neq\ga_t$.
\end{enumerate}
In words, an open path may jump using arrows and cannot stay at the same site
if there is a blocking symbol at such a site. For $i,j\in\La$ and $s\leq u$,
we write $(i,s)\leadsto(j,u)$ if there is an open path $\ga$ that leads from
$(i,s)$ to $(j,u)$, i.e., $\ga_s=i$ and $\ga_u=j$. The following lemma, whose
formal proof will be given in Section~\ref{S:percproof}, gives the promised
percolation representation of additive Markov processes in $S=\Pc_{\rm
  dec}(\La)$, and their duals. Moreover, formula (\ref{duint}) gives a
graphical interpretation of the pathwise duality of Theorem~\ref{T:addual}.
Note that in the present setting, the duality function takes the form
$\psi(x,y)=1_{\{x\sub y^{\rm c}\}}=1_{\{x\cap y=\emptyset\}}$.

\bl[Percolation representation]\label{L:percrep}
Almost surely, for all $s\leq t$ and $x\in S$,
\be\label{percrep}
\Xb_{s,u}(x)=\{j\in\La:(i,s)\leadsto(j,u)\mbox{ for some }i\in x\},
\ee
and the left-continuous version of the dual stochastic flow is a.s.\ given by
\be\label{dualrep}
\Yb_{s-,u-}(y)=\{j\in\La:(j,-u)\leadsto(i,-s)\mbox{ for some }i\in y\}
\qquad(s\leq u,\ y\in S').
\ee
Moreover, for deterministic times $s\leq u$, a.s.
\be\label{duint}
1_{\txt\{\Xb_{s,t-}(x)\cap\Yb_{-u,-t}(y)=\emptyset\}}
=1_{\txt\{(i,s)\not\leadsto(j,u)\ \forall\ i\in x,\ j\in y\}}
\qquad(s\leq t\leq u).
\ee
\el

By Birkhoff's representation theorem, Lemmas~\ref{L:madd}
and\ref{L:percrep} apply to additive processes taking values in general
distributive lattices. For nondistributive lattices, the picture is not as
nice, but we can still show that each additive process has a percolation
representation. The next lemma will be proved in Section~\ref{S:percproof}.

\bl[Additive maps on general lattices]\label{L:nondis}
Let $S$ be a finite lattice. Then $S$ is $(0,\vee)$-isomorphic to a
join-semilattice of sets, i.e., there exists a finite set $\La$ and a set
$T\sub\Pc(\La)$ such that $\emptyset\in T$ and $T$ is closed under unions, and
$S\cong T$. Moreover, if $T$ is such a join-semilattice of sets, then each
additive map $m:T\to T$ can be extended to an additive map $\ov
m:\Pc(\La)\to\Pc(\La)$.
\el

By Lemma~\ref{L:nondis}, if $G$ is the generator of a Markov process $X$
in $T$, and $X$ is additive in the sense that $G$ has a random mapping
representation of the form (\ref{GPois}) where all maps $m\in\Gi$ are
additive, then we can extend $X$ to a Markov process in $\Pc(\La)$
with generator of the form
\be\label{Gext}
\ov Gf(x)=\sum_{m\in\Gi}r_m\big(f(\ov m(x))-f(x)\big)\qquad\big(x\in\Pc(\La)\big).
\ee
This extended process leaves $T$ invariant, and by Lemma~\ref{L:nondis}, it has
a percolation representation of the form (\ref{percrep}). If $S$ is
nondistributive, however, then $T$ can never be chosen such that it is also
closed under intersections, which means that $\{x^{\rm c}:x\in T\}$ is not
$(0,\vee)$-isomorphic to $S'$, and the parts of Lemma~\ref{L:percrep}
referring to the dual process break down. Of course, the extended process
(\ref{Gext}) has a dual that can be interpreted in terms of open paths, but it
is unclear if the extended process can ever be chosen in such a way that it
leaves a subspace invariant that is $(0,\vee)$-isomorphic to $S'$.

\subsection{Some open problems}\label{S:open}

The problem how to decide whether a given Markov generator can be represented
in monotone or additive maps has already been mentioned in the text. These
problems have been open for a long time and appear to be hard.

We have also only partially resolved the question whether additive systems
taking values in a nondistributive lattice have a ``nice'' percolation
representation together with their dual. If $S$ is a nondistributive lattice,
then by Lemma~\ref{L:nondis}, $S$ is $(0,\vee)$-isomorphic to a
join-semilattice of sets, but by Birkhoff's theorem such a join-semilattice of
sets can never be closed under intersections. This also means that we cannot
simultaneously represent the dual lattice $S'$ as a join-semilattice of sets
on the same space, so that the duality map $x\mapsto x'$ is the complement map.
However, to get a useful percolation representation of $S$ and $S'$ together,
it would suffice to represent $S$ and $S'$ as join-semilattice of sets on the
same space in such a way that the duality map $\li x,y\re$ takes the form
$1_{\{x\cap y=\emptyset\}}$. It is not clear if this can be done.

A more urgent question is perhaps if all this abstract theory ``is actually
good for anything''.  Additive systems duality is clearly a very useful
tool, and as the two-stage contact process shows (see Section~\ref{S:Krone}
below), one sometimes needs more general lattices than the Boolean algebra
$\{0,1\}^\La$. Our general approach to constructing percolation
  representations for such lattices seems to be new. The more general but
also more complicated duality for monotone systems, that are not necessarily
additive, has so far found few applications, although Gray \cite{Gra86} did
use it to prove nontrivial statements and we used it in our previous work
  \cite{SS14} to derive a useful subduality for systems with cooperative
  branching.

For many additive systems including the contact process, using duality, it is
fairly easy to prove that starting from an arbitrary translation invariant
initial law, the law of the system converges to a convex combination of the
upper invariant law and the delta measure on the zero configuration; see
\cite[Thm~III.5.18]{Lig85}. As far as we know, it is an open problem to
generalize these techniques to monotone systems that are not additive, such as
the systems with cooperative branching discussed in Section~\ref{S:coop}
below. We hope that our present systematic treatment of monotone systems
duality can contribute to such an undertaking.

In this context, we mention one more open problem. In Section~\ref{S:monpart}
below, we generalize Theorem~\ref{T:mondual} about monotone systems duality to
infinite underlying spaces. We only show, however, that the dual process
is well-defined started from finite initial states. It is an open problem to
construct the dual process with infinite initial states, as would be needed,
e.g., to study invariant laws of the dual process.

\section{Examples}\label{S:examples}

\subsection{Siegmund's duality}\label{S:Siegmund}

In this subsection, we show that for finite state spaces, Siegmund's duality
\cite{Sie76} is a special case of Theorem~\ref{T:addual} which gives a
pathwise dual for additive Markov processes taking values in a finite lattice.

Let $S=\{0,\ldots,n\}$ $(n\geq 2)$ be a finite totally ordered set with at
least two elements, and let $S'$ denote the same space, equipped with the
reversed order, which is a dual of $S$ in the sense defined in
Section~\ref{S:addual}. Since $x\vee y$ equals either $x$ or $y$ for each
$x,y\in S$, it is easy to see that a map $m:S\to S$ is additive if and only if
$m$ is monotone and $m(0)=0$. By Lemma~\ref{L:addual}, for each such map there
exists a unique map $m':S\to S$ that is dual to $m$ with respect to the
duality function $\psi(x,y)=1_{\{x\leq y\}}$. This map $m'$ is additive viewed
as a map on $S'$, i.e., $m'$ is monotone and $m'(n)=n$. Theorem~\ref{T:addual}
now gives us the following result. We note that Siegmund's result \cite{Sie76}
applies more generally to processes taking values in a real interval, which
includes, e.g., a duality between Brownian motions on $\half$ with absorption
at 0 and those with reflection at 0.

\bt[Pathwise Siegmund's duality]\label{T:Siegmund}
Let $S=\{0,\ldots,n\}$ $(n\geq 2)$ be a finite totally ordered set and let $X$
be a Markov process in $S$ whose generator has a random mapping representation
of the form
\be
Gf(x)=\sum_{m\in\Gi}r_m\big(f(m(x))-f(x)\big)\qquad(x\in S),
\ee
where all maps $m\in\Gi$ are monotone and satisfy $m(0)=0$. Then the Markov
process $Y$ in $S'$ with generator
\be
Hf(y):=\sum_{m\in\Gi}r_m\big(f(m'(y))-f(y)\big)\qquad(y\in S')
\ee
is pathwise dual to $X$ with respect to the duality function
$\psi(x,y):=1_{\{x\leq y\}}$.
\et

A Markov process $X$ taking values in a partially ordered set $S$ is
called \emph{monotone} if its transition kernels $P_t$ are monotone for each
$t\geq 0$. If $S=\{0,\ldots,n\}$ is totally ordered, then by
Proposition~\ref{P:monsuf}, it follows that $P_t$ is monotonically representable for each $t\geq 0$, and by Lemma~\ref{L:genrep} the same is true for
the generator $G$. If $0$ is a trap (i.e., the rate of jumps away from $0$
is zero), then each random mapping representation for $G$ involves only maps
satisfying $m(0)=0$. Since each pathwise dual is also a dual in the classical
sense (\ref{dual}), Theorem~\ref{T:Siegmund} allows us to conclude:

\bp[Siegmund's duality]\label{P:Siegmund}
Let $S=\{0,\ldots,n\}$ $(n\geq 2)$ be a finite totally ordered set and let $X$
be a monotone Markov process in $S$ for which $0$ is a trap. Then there exists
a monotone Markov process $Y$ in $S$ for which $n$ is a trap, such that
\be\label{Sied}
\P[X_t\leq Y_0]=\P[X_0\leq Y_t]\qquad(t\geq 0)
\ee
for arbitrary deterministic initial states $X_0$ and $Y_0$.
\ep

It is easy to see that (\ref{Sied}) determines the transition probabilities of
$Y$ uniquely, even though there usually is more than one way to represent the
generators of $X$ and $Y$ in terms of monotone maps and hence to construct a
coupling that realizes the pathwise duality of Theorem~\ref{T:Siegmund}. The
fact that Siegmund's duality can always be realized as a pathwise duality has
been observed before in \cite{CS85}, and was in fact the main point of that
paper.

%For example, a random walk can be constructed by Poisson events that say that
%the process should jump up by one step, regardless of the present state, but
%we can also construct the same process by Poisson events that apply only to
%one transition between two given states.

\subsection{Additive interacting particle systems}\label{S:addpart}

Additive particle systems in the classical sense of \cite{Gri79,Har78} are
additively representable Markov processes $X$ with state space of the form
$S=\Pc(\La)$ where $\La$ is a countable set. Concentrating on finite $\La$ for
the moment, this is a special case of the set-up from Section~\ref{S:percol},
where the set $\La$ from Lemmas~\ref{L:madd} and \ref{L:percrep} is equipped
with the trivial order, so that $\Pc_{\rm dec}(\La)=\Pc(\La)$, and the
conditions (\ref{Mprop}) are void. In this case, (\ref{minM}) simply defines a
one-to-one correspondence between additive maps $m:S\to S$ and sets
$M\sub\La\times\La$, and Lemma~\ref{L:percrep} shows that $X$ and its dual $Y$
can be represented in terms of the same percolation substructure and are
pathwise dual in the sense of (\ref{duint}).

The well-known voter model, contact process, and exclusion process are
additive particle systems that fall into this framework. Here, we illustrate
our framework on the concrete example of the voter model. We refer to
\cite{Swa16} for an analogue treatment of contact and exclusion processes
which can be represented in terms of birth maps, death maps, and exclusion
maps that are all additive. A \emph{voter model} is a Markov process that
has a generator of the form
\be\label{Gvoter}
G_{\rm voter}f(x):=\sum_{i,j\in\La}r_{ij}
\big(f\big({\rm vot}_{ij}(x)\big)-f\big(x\big)\big)
\qquad\big(x\in\Pc(\La)\big),
\ee
where we define \emph{voter model maps}
${\rm vot}_{ij}:\Pc(\La)\to\Pc(\La)$ by
\be\label{votmap}
{\rm vot}_{ij}(x):=\left\{\ba{ll}
\dis x\cup\{j\}\quad&\dis\mbox{if }i\in x,\\[5pt]
\dis x\beh\{j\}\quad&\dis\mbox{if }i\not\in x,\ea\right.
\qquad\big(x\in\Pc(\La),\ i,j\in\La\big).
\ee
The map ${\rm vot}_{ij}$ corresponds (in the sense of (\ref{minM})) to the set
$M_{ij}\sub\La\times\La$ given by $M_{ij}=\{(k,k):k\in\La,\ k\neq
j\}\cup\{(i,j)\}$, which is represented by an arrow from $i$ to $j$ and
simultaneously a blocking symbol at $j$.

By Lemma~\ref{L:madd}, the dual map ${\rm vot}'_{ij}$ corresponds to the set
$M'=\{(k,k):k\in\La,\ k\neq j\}\cup\{(j,i)\}$, which through (\ref{minM}),
$M'$ defines the dual map
\be\label{rwmap}
{\rm rw}_{ji}(y):=\left\{\ba{ll}
\dis(y\beh\{j\}\big)\cup\{i\}\quad&\dis\mbox{if }j\in y,\\[5pt]
\dis y\quad&\dis\mbox{if }j\not\in y,\ea\right.
\qquad\big(y\in\Pc(\La),\ i,j\in\La\big),
\ee
which is represented by an arrow from $j$ to $i$ and simultaneously a
blocking symbol at $j$. We interpret ${\rm rw}_{ji}$ as a \emph{coalescing
  random walk map}, i.e., when ${\rm rw}_{ji}$ is applied, if there is a
particle at $j$, then this particle jumps to $i$, coalescing with any particle
that may already be present.
 
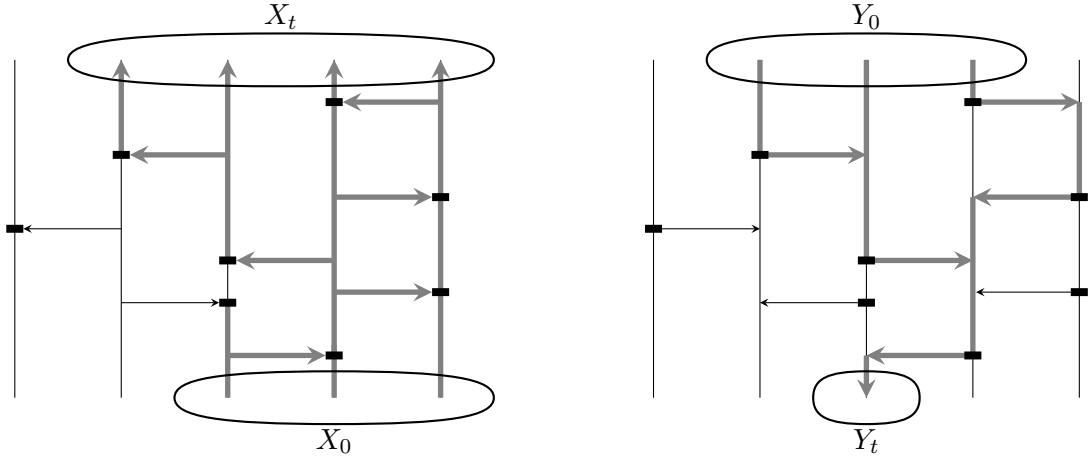
\begin{figure}[tcb]
\begin{center}
\begin{tikzpicture}[scale=1.4,>=stealth]
\begin{scope}
\foreach \x in {1,...,5}
\draw (\x,-0.1)--(\x,3.1);

\draw[->,line width=2pt,gray] (3,0.3)--++(0.92,0);
\draw[->] (2,0.8)--++(0.92,0);
\draw[->,line width=2pt,gray] (4,0.9)--++(0.92,0);
\draw[->,line width=2pt,gray] (4,1.2)--++(-0.92,0);
\draw[->] (2,1.5)--++(-0.92,0);
\draw[->,line width=2pt,gray] (4,1.8)--++(0.92,0);
\draw[->,line width=2pt,gray] (3,2.2)--++(-0.92,0);
\draw[->,line width=2pt,gray] (5,2.7)--++(-0.92,0);

\draw[line width=2pt,gray,->] (2,2.2)--(2,3.1);
\draw[line width=2pt,gray] (4,-0.1)--(4,0.4);
\draw[line width=2pt,gray] (3,-0.1)--(3,0.8);
\draw[line width=2pt,gray,->] (3,1.2)--(3,3.1);
\draw[line width=2pt,gray] (4,0.3)--(4,2.7);
\draw[line width=2pt,gray,->] (4,2.7)--(4,3.1);
\draw[line width=2pt,gray] (5,-0.1)--(5,0.9);
\draw[line width=2pt,gray] (5,0.9)--(5,1.8);
\draw[line width=2pt,gray,->] (5,1.8)--(5,3.1);

\draw[line width=3pt] (4,0.3)++(-0.08,0)--++(0.16,0);
\draw[line width=3pt] (3,0.8)++(-0.08,0)--++(0.16,0);
\draw[line width=3pt] (5,0.9)++(-0.08,0)--++(0.16,0);
\draw[line width=3pt] (3,1.2)++(-0.08,0)--++(0.16,0);
\draw[line width=3pt] (1,1.5)++(-0.08,0)--++(0.16,0);
\draw[line width=3pt] (5,1.8)++(-0.08,0)--++(0.16,0);
\draw[line width=3pt] (2,2.2)++(-0.08,0)--++(0.16,0);
\draw[line width=3pt] (4,2.7)++(-0.08,0)--++(0.16,0);

\draw[thick] plot[smooth cycle,tension=1.3]
coordinates{(1.5,3.1) (3.5,2.85) (5.5,3.1) (3.5,3.35)};
\draw[thick] plot[smooth cycle,tension=1.3]
coordinates{(2.5,-0.1) (4,-0.35) (5.5,-0.1) (4,0.15)};

\draw (3.5,3.3) node[black,above] {$X_t$};
\draw (4,-0.3) node[black,below] {$X_0$};
\end{scope}

\begin{scope}[xshift=6cm]
\foreach \x in {1,...,5}
\draw (\x,-0.1)--(\x,3.1);

\draw[->,line width=2pt,gray] (4,0.3)--++(-1,0);
\draw[->] (3,0.8)--++(-1,0);
\draw[->] (5,0.9)--++(-0.97,0);
\draw[->,line width=2pt,gray] (3,1.2)--++(1,0);
\draw[->] (1,1.5)--++(1,0);
\draw[->,line width=2pt,gray] (5,1.8)--++(-1,0);
\draw[->,line width=2pt,gray] (2,2.2)--++(1,0);
\draw[->,line width=2pt,gray] (4,2.7)--++(1,0);

\draw[line width=2pt,gray] (2,2.2)--(2,3.1);
\draw[<-,line width=2pt,gray] (3,-0.1)--(3,0.3);
\draw[line width=2pt,gray] (3,1.2)--(3,3.1);
\draw[line width=2pt,gray] (4,0.3)--(4,1.8);
\draw[line width=2pt,gray] (4,2.7)--(4,3.1);
\draw[line width=2pt,gray] (5,1.8)--(5,2.7);

\draw[line width=3pt] (4,0.3)++(-0.08,0)--++(0.16,0);
\draw[line width=3pt] (3,0.8)++(-0.08,0)--++(0.16,0);
\draw[line width=3pt] (5,0.9)++(-0.08,0)--++(0.16,0);
\draw[line width=3pt] (3,1.2)++(-0.08,0)--++(0.16,0);
\draw[line width=3pt] (1,1.5)++(-0.08,0)--++(0.16,0);
\draw[line width=3pt] (5,1.8)++(-0.08,0)--++(0.16,0);
\draw[line width=3pt] (2,2.2)++(-0.08,0)--++(0.16,0);
\draw[line width=3pt] (4,2.7)++(-0.08,0)--++(0.16,0);

\draw[thick] plot[smooth cycle,tension=1.3]
coordinates{(1.5,3.1) (3,2.85) (4.5,3.1) (3,3.35)};
\draw[thick] plot[smooth cycle,tension=1.3]
coordinates{(2.5,-0.1) (3,-0.35) (3.5,-0.1) (3,0.15)};

\draw (3,3.3) node[black,above] {$Y_0$};
\draw (3,-0.3) node[black,below] {$Y_t$};
\end{scope}
\end{tikzpicture}
\caption{Graphical representation of a voter model and its dual.}
\label{fig:voter}
\end{center}
\end{figure}

By Theorem~\ref{T:addual}, the voter model $X$ is pathwise dual to the system
of coalescing random walks $Y$ with generator
\be\label{Grw}
G_{\rm rw}f(y):=\sum_{i,j\in\La}r_{ij}
\big(f\big({\rm rw}_{ji}(y)\big)-f\big(y\big)\big)
\qquad\big(y\in\Pc(\La)\big),
\ee
and by Lemma~\ref{L:percrep}, the stochastic flows associated with $X$ and $Y$
can be represented in terms of open paths in a percolation substructure.
In Figure~\ref{fig:voter}, we have drawn an example of such a percolation
substructure, together with a voter model (in the upward time direction) and
system of coalescing random walks (with time running downwards).
In the picture for the dual process $Y$, we have reversed the direction of all
Poisson arrows, in line with formula (\ref{Madj}) for the dual of an additive
map.

\subsection{Krone's duality}\label{S:Krone}

Steve Krone \cite{Kro99} has studied a two-stage contact
process, which is a Markov process with state space of the form
$S=\{0,1,2\}^\La$. The main interest is in the case $\La=\Z^d$ but the
construction works in the same way if $\La$ is finite. If $x(i)=0,1$, or $2$,
then he interprets this as the site $i$ being occupied by no individual, a
young individual, or an adult individual, respectively. For each $i,j\in\La$,
consider the maps $a_i,b_{ij},c_i,d_i,e_i$ defined by
\be\ba{lr@{\,}c@{\,}ll}\label{Kromaps}
\mbox{grow up}\quad&a_i(x)(k)&:=&2\quad\mbox{if }k=i,\ x(i)=1,
\quad&:=x(k)\mbox{ otherwise,}\\[5pt] 
\mbox{give birth}\quad&b_{ij}(x)(k)&:=&1\quad\mbox{if }k=j,\ x(i)=2,\ x(j)=0,
\quad&:=x(k)\mbox{ otherwise,}\\[5pt] 
\mbox{young dies}\quad&c_i(x)(k)&:=&0\quad\mbox{if }k=i,\ x(i)=1,
\quad&:=x(k)\mbox{ otherwise,}\\[5pt] 
\mbox{death}\quad&d_i(x)(k)&:=&0\quad\mbox{if }k=i,
\quad&:=x(k)\mbox{ otherwise,}\\[5pt] 
\mbox{grow younger}\quad&e_i(x)(k)&:=&1\quad\mbox{if }k=i,\ x(i)=2,
\quad&:=x(k)\mbox{ otherwise.}
\ec
Except for the last one, all these maps have natural biological
interpretations. For example, the map $b_{ij}$ has the effect that if $i$ is
occupied by an adult individual and $j$ is empty, then the adult individual at
$i$ gives birth to a young individual at $j$. Applying the map $c_i$ with an
appropriate rate models the commonly observed fact that  young individuals  die
at a higher rate than adults.

We set $S':=S$ and define a map $S\ni x\mapsto x'\in S'$ by
\be\label{Kroac}
x'(i):=2-x(i)\qquad(i\in\La).
\ee
Then $S'$, together with the map $x\mapsto x'$, is a dual of $S$ in the sense
of (\ref{duspa}). The set $S$, being the product of totally ordered sets, is a
distributive lattice and so we can apply Theorem~\ref{T:addual} to find pathwise
duals of additive Markov processes in $S$. Our choice of the dual space $S'$
means that the duality function $\psi(x,y)=\li x,y\re$ from (\ref{tetSS})
takes the form
\be\label{Kropsi}
\psi(x,y)=1_{\txt\{x\leq y'\}}=1_{\txt\{x(i)+y(i)\leq 2\ \forall\ i\in\La\}}.
\ee
The good news is that the maps in (\ref{Kromaps}) are all additive.

\bl[Additive maps]\label{L:Krodu}
The maps in (\ref{Kromaps}) are all additive and their dual maps with respect
to the duality function in (\ref{Kropsi}) are given by
\be\label{Krodu}
a'_i=a_i,\quad b'_{ij}=b_{ji},\quad c'_i=e_i,\quad d'_i=d_i,\quad e'_i=c_i.
\ee
\el

Lemma~\ref{L:Krodu} can be checked by straightforward, but somewhat lengthy
considerations. Things become easier if we construct a percolation
representation for two-stage contact processes. We equip $\La$ with the
trivial order, we view $\{0,1\}$ 
as a totally ordered set with two elements,
and we equip $\La\times\{0,1\}$ with the product order. Then
\be
S\cong\Pc_{\rm dec}\big(\La\times\{0,1\}\big)
\quand
S'\cong\Pc_{\rm inc}\big(\La\times\{0,1\}\big),
\ee
where for the duality map we now take the complement map $x\mapsto x':=x^{\rm
  c}$ so that in this new representation the duality function takes the usual
form $\psi(x,y)=1_{\{x\cap y=\emptyset\}}$. To compare this with our earlier
representation of the process, note that for $x\in\Pc_{\rm
  dec}(\La\times\{0,1\})$ and $i\in\La$,
\be\label{forwrep}
\{\sig:(i,\sig)\in x\}=\emptyset,\quad\{0\},\quad\mbox{or}\quad\{0,1\},
\ee
which we interpret as $x(i)=0,1$, or $2$ in our previous representation of $S$.
Likewise, for $y\in\Pc_{\rm inc}(\La\times\{0,1\})$ and $i\in\La$,
\be\label{backrep}
\{\sig:(i,\sig)\in y\}=\emptyset,\quad\{1\},\quad\mbox{or}\quad\{0,1\},
\ee
which we interpret as $y(i)=0,1$, or $2$ in our previous representation of~$S'$.

\begin{figure}[tcb]
\begin{center}
\begin{tikzpicture}[scale=1.7,>=stealth]
\begin{scope}[yscale=0.75]
\foreach \x in {1,...,4}
 \draw (\x,0) --++(0,4);
\foreach \x in {1,...,4}
 \draw (\x,0) ++(0.3,0) --++(0,4);

\coordinate (a) at (2,0.3);
\coordinate (b) at (3,0.9);
\coordinate (c) at (4,0.9);
\coordinate (d) at (2,1.2);
\coordinate (e) at (2,1.8);
\coordinate (f) at (3,1.8);
\coordinate (g) at (2,2.5);
\coordinate (h) at (3,2.7);
\coordinate (i) at (2,3.0);
\coordinate (j) at (3,3.1);
\coordinate (k) at (4,3.1);
\coordinate (l) at (1,2.5);
\coordinate (m) at (1,3.5);
\coordinate (n) at (2,3.6);
\coordinate (o) at (3,3.6);
\coordinate (p) at (4,3.7);
\coordinate (q) at (4,1.7);

\draw[line width=2pt,gray] (2,0)--(i);
\draw[line width=2pt,gray] (a)++(0.3,0)--($(i)+(0.3,0)$);
\draw[line width=2pt,gray] (l)--(m);
\draw[line width=2pt,gray,->] (n)--(2,4);
\draw[line width=2pt,gray,->] (b)--(3,4);
\draw[line width=2pt,gray,->] (h)++(0.3,0)--(3.3,4);
\draw[line width=2pt,gray] (k)--(p);
\draw[line width=2pt,gray] (4,0)--(q);
\draw[line width=2pt,gray] (4.3,0)--($(q)+(0.3,0)$);

\draw[line width=3pt] (d)++(-0.08,0)--++(0.16,0);
\draw[line width=3pt] (i)++(-0.08,0)--++(0.16,0);
\draw[line width=3pt] (i)++(0.3,0)++(-0.08,0)--++(0.16,0);
\draw[line width=3pt] (m)++(-0.08,0)--++(0.16,0);
\draw[line width=3pt] (m)++(0.3,0)++(-0.08,0)--++(0.16,0);
\draw[line width=3pt] (p)++(-0.08,0)--++(0.16,0);
\draw[line width=3pt] (q)++(-0.08,0)--++(0.16,0);
\draw[line width=3pt] (q)++(0.3,0)++(-0.08,0)--++(0.16,0);

\draw[<->,line width=2pt,gray] (a)--++(0.3,0);
\draw[->,line width=2pt,gray] (c)++(0.3,0)--(b);
\draw[->,line width=2pt,gray] (d)++(0.3,0)--++(-0.22,0);
\draw[->,line width=2pt,gray] (e)++(0.3,0)--++(0.7,0);
\draw[->,line width=2pt,gray] (g)++(0.3,0)--(l);
\draw[<->,line width=2pt,gray] (h)--++(0.3,0);
\draw[->,line width=2pt,gray] (j)++(0.3,0)--++(0.7,0);
\draw[->,line width=2pt,gray] (o)++(0.3,0)--(n);
\draw[->,thick] (p)++(0.3,0)--++(-0.22,0);

\draw (1.15,0) node[below] {0};
\draw (2.15,0) node[below] {1};
\draw (3.15,0) node[below] {0};
\draw (4.15,0) node[below] {2};
\draw (1.15,4) node[above] {0};
\draw (2.15,4) node[above] {1};
\draw (3.15,4) node[above] {2};
\draw (4.15,4) node[above] {0};

\draw (2.65,-0.3) node[below] {$X_0$};
\draw (2.65,4.3) node[above] {$X_t$};

%\node[green] at (2,0.3) {a};
%\node[green] at (3,0.9) {b};
%\node[green] at (4,0.9) {c};
%\node[green] at (2,1.2) {d};
%\node[green] at (2,1.8) {e};
%\node[green] at (3,1.8) {f};
%\node[green] at (2,2.5) {g};
%\node[green] at (3,2.7) {h};
%\node[green] at (2,3.0) {i};
%\node[green] at (3,3.1) {j};
%\node[green] at (4,3.1) {k};
%\node[green] at (1,2.5) {l};
%\node[green] at (1,3.5) {m};
%\node[green] at (2,3.6) {n};
%\node[green] at (3,3.6) {o};
%\node[green] at (4,3.7) {p};
%\node[green] at (4,1.7) {q};
\end{scope}

\begin{scope}[xshift=5cm,yscale=0.75]
\foreach \x in {1,...,4}
 \draw (\x,0) --++(0,4);
\foreach \x in {1,...,4}
 \draw (\x,0) ++(0.3,0) --++(0,4);

\coordinate (a) at (2,0.3);
\coordinate (b) at (3,0.9);
\coordinate (c) at (4,0.9);
\coordinate (d) at (2,1.2);
\coordinate (e) at (2,1.8);
\coordinate (f) at (3,1.8);
\coordinate (g) at (2,2.5);
\coordinate (h) at (3,2.7);
\coordinate (i) at (2,3.0);
\coordinate (j) at (3,3.1);
\coordinate (k) at (4,3.1);
\coordinate (l) at (1,2.5);
\coordinate (m) at (1,3.5);
\coordinate (n) at (2,3.6);
\coordinate (o) at (3,3.6);
\coordinate (p) at (4,3.7);
\coordinate (q) at (4,1.7);

\draw[line width=2pt,gray] (1,4)--(m);
\draw[line width=2pt,gray] (1.3,4)--($(m)+(0.3,0)$);
\draw[line width=2pt,gray,->] (a)--(2,0);
\draw[line width=2pt,gray,->] ($(e)+(0.3,0)$)--(2.3,0);
\draw[line width=2pt,gray,->] (h)--(3,0);
\draw[line width=2pt,gray,->] (3.3,4)--(3.3,0);
\draw[line width=2pt,gray] (4,4)--(p);
\draw[line width=2pt,gray] (4.3,4)--($(q)+(0.3,0)$);
\draw[line width=2pt,gray,->] ($(c)+(0.3,0)$)--(4.3,0);

\draw[line width=3pt] (d)++(-0.08,0)--++(0.16,0);
\draw[line width=3pt] (i)++(-0.08,0)--++(0.16,0);
\draw[line width=3pt] (i)++(0.3,0)++(-0.08,0)--++(0.16,0);
\draw[line width=3pt] (m)++(-0.08,0)--++(0.16,0);
\draw[line width=3pt] (m)++(0.3,0)++(-0.08,0)--++(0.16,0);
\draw[line width=3pt] (p)++(-0.08,0)--++(0.16,0);
\draw[line width=3pt] (q)++(-0.08,0)--++(0.16,0);
\draw[line width=3pt] (q)++(0.3,0)++(-0.08,0)--++(0.16,0);

\draw[<->,line width=2pt,gray] (a)--++(0.3,0);
\draw[<-,line width=2pt,gray] (c)++(0.3,0)--(b);
\draw[<-,thick] (d)++(0.3,0)--++(-0.22,0);
\draw[<-,line width=2pt,gray] (e)++(0.3,0)--++(0.7,0);
\draw[<-,thick] (g)++(0.3,0)--(l);
\draw[<->,line width=2pt,gray] (h)--++(0.3,0);
\draw[<-,thick] (j)++(0.3,0)--++(0.7,0);
\draw[<-,thick] (o)++(0.3,0)--(n);
\draw[<-,line width=2pt,gray] (p)++(0.3,0)--++(-0.22,0);

\draw (1.15,0) node[below] {0};
\draw (2.15,0) node[below] {2};
\draw (3.15,0) node[below] {2};
\draw (4.15,0) node[below] {1};
\draw (1.15,4) node[above] {2};
\draw (2.15,4) node[above] {0};
\draw (3.15,4) node[above] {1};
\draw (4.15,4) node[above] {2};

\draw (2.65,-0.3) node[below] {$Y_t$};
\draw (2.65,4.3) node[above] {$Y_0$};

%\node[green] at (2,0.3) {a};
%\node[green] at (3,0.9) {b};
%\node[green] at (4,0.9) {c};
%\node[green] at (2,1.2) {d};
%\node[green] at (2,1.8) {e};
%\node[green] at (3,1.8) {f};
%\node[green] at (2,2.5) {g};
%\node[green] at (3,2.7) {h};
%\node[green] at (2,3.0) {i};
%\node[green] at (3,3.1) {j};
%\node[green] at (4,3.1) {k};
%\node[green] at (1,2.5) {l};
%\node[green] at (1,3.5) {m};
%\node[green] at (2,3.6) {n};
%\node[green] at (3,3.6) {o};
%\node[green] at (4,3.7) {p};
%\node[green] at (4,1.7) {q};
\end{scope}
\end{tikzpicture}
\caption{Graphical representation of a two-stage contact process and its dual.}
\label{fig:Krone}
\end{center}
\end{figure}
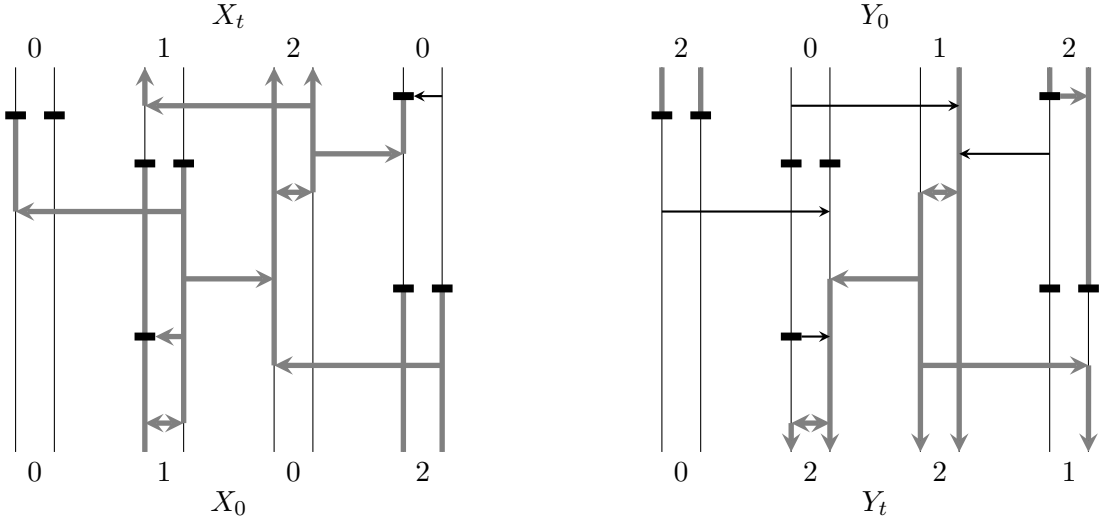

An example of a percolation representation for a two-stage contact process and
its dual are shown in Figure~\ref{fig:Krone}. Here, the maps from
(\ref{Kromaps}) are represented as follows:
\begin{itemize}
\item[$a_i$] arrows between $(i,0)$ and $(i,1)$ (in both directions conforming
  to the rules of Lemma~\ref{L:madd}),
\item[$b_{ij}$] arrow from $(i,1)$ to $(j,0)$,
\item[$c_i$] blocking symbol at $(i,0)$, arrow from $(i,1)$ to $(i,0)$,
\item[$d_i$] blocking symbols at $(i,0)$ and $(i,1)$,
\item[$e_i$] blocking symbol at $(i,1)$, arrow from $(i,1)$ to $(i,0)$.
\end{itemize}

\subsection{Gray's duality}\label{S:Gray}

Gray \cite{Gra86} proved a general duality for monotone spin systems that
need, in general, not be additive. Although Gray's formulation of the duality
differs considerably from ours, in the present subsection, we will show that
his dual is a special case of the process with generator $H_\dott$, which is
defined as in Theorem~\ref{T:mondual} but with the dual maps $m^\ast$ replaced
by the analogous maps $m^\dott$ from Lemma~\ref{L:altmaps}.

Gray considers spin systems, which are Markov processes with
state space $S:=\{0,1\}^\La$ and generator of the form
\be\label{GGray}
Gf(x)=\sum_{i\in\La}\bet_i(x)\big(f(x\vee\eps_i)-f(x)\big)
+\sum_{i\in\La}\de_i(x)\big(f(x\wedge(1-\eps_i))-f(x)\big),
\ee
where $\eps_i(j):=1_{\{i=j\}}$. Gray's emphasis is on the case $\La=\Z$ but the
arguments are the same for finite $\La$. Gray assumes that his systems are
\emph{attractive}, which means that the function $\bet_i:S\to\half$ is
monotone and the function $\de_i:S\to\half$ is \emph{anti-monotone}, i.e.,
$-\de_i$ is monotone. He then proves the following fact.

\bl[Attractive spin systems are  monotonically representable]\label{L:spinmon}
Let $G$ be of the form (\ref{GGray}) and assume that for each $i\in\La$, the
function $\bet_i$ is monotone and $\de_i$ is anti-monotone. Then there exists a
set $\Gi$ of monotone maps $m:S\to S$ and nonnegative
constants $(r_m)_{m\in\Gi}$ such that $G$ has the form
\be\label{GPois2}
Gf(x)=\sum_{m\in\Gi}r_m\big(f(m(x))-f(x)\big)\qquad(x\in S).
\ee
\el

Using the random mapping representation (\ref{GPois2}) of $G$ in terms of
monotone maps, Gray then constructs a Poisson point set $\De$ whose elements
are pairs $(m,t)$ with $m\in\Gi$ and $t\in\R$, and uses this to define a
stochastic flow $(\Xb_{s,t})_{s\leq t}$ as in (\ref{Phidef}).\footnote{See
  formulas (10)--(12) of \cite{Gra86}. His notation for $\Xb_{s,t}(x)$ is
  $\xi(s,t,x)$.}  Gray then defines an \emph{$[s,u]$-path} from $x$ to $y$
(with $x,y\in S$) to be a cadlag function $\pi:[s,u]\to S$ with $\pi(s)=x$ and
$\pi(u)=y$, such that:\footnote{Gray does not state his definition entirely
  correctly. In his definition of minimality, he also includes, probably
  inadvertently, the condition that $\ti\pi(s)=x$.}
\begin{enumerate}
\item $\Xb_{t,t'}(\pi(t))\geq\pi(t')$ for all $s\leq t\leq t'\leq u$.
\item $\pi$ is \emph{minimal} in the sense that if a cadlag function
  $\ti\pi:[s,u]\to S$ with $\ti\pi(u)=y$ satisfies (i) and
  $\ti\pi(t)\leq\pi(t)$ for all $t\in[s,u]$, then $\pi=\ti\pi$.
\end{enumerate}
Next, for each $s\leq u$, he defines random maps from $S$ to $\Pc(S)$ by
\be\label{zetadef}
\zeta_{s,u}(y)
:=\big\{x\in S:\mbox{ there exists an $[s,u]$-path from $x$ to $y$}\big\}
\qquad(y\in S).
\ee
With this notation, Theorem~1 in \cite{Gra86} reads:

\bt[Gray's duality]\label{T:Gray}
For each $x,y\in S$ and $s\leq u$, one has $y\leq\Xb_{s,u}(x)$ if and only if
$z\leq x$ for some $z\in\zeta_{s,u}(y)$.
\et

To relate this to our work, for each monotone map $m:S\to S$, let $m^\dott$
denote the dual map defined in (\ref{mstar}), where for $S'$ we choose the set
$S$ equipped with the reversed order, i.e.,
\be\label{mstar2}
m^\dott(B):=\bigcup_{x\in B}(m^{-1}(\{x\}^\up))_{\rm min}
\qquad\big(B\in\Pc(S)\big),
\ee
where on the right-hand side $\{x\}^\up$ and the minimum are taken with
respect to the order on $S$. Set $\De^\dott:=\{(m^\dott,-t):(m,t)\in\De\}$ and
use this Poisson set to define a stochastic flow $(\Yb^\dott_{s,t})_{s\leq t}$
that by Proposition~\ref{P:pathdual} and Lemma~\ref{L:altmaps} is dual to
$(\Xb_{s,t})_{s\leq t}$ with respect to the duality function $\ti\phi$ from
(\ref{altphi}). Then Theorem~\ref{T:Gray} is an immediate consequence of this
duality and the following fact, that will be proved in
Section~\ref{S:Grayproof}.

\bp[Reformulation of Gray's dual]\label{P:Graystar}
For each $y\in S$ and $s\leq u$, one has
$\zeta_{s,u}(y)=\Yb^\dott_{-u,-s}\big(\{y\}\big)$.
\ep

Note that the duality of $(\Yb^\dott_{s,t})_{s\leq t}$ to
$(\Xb_{s,t})_{s\leq t}$ with respect to the duality function from
(\ref{altphi}) implies that for each $x,y\in S$ and $s\leq  t$,
\be
1_{\txt\{x\geq z\mbox{ for some }z\in\Yb^\dott_{-t,-s}(\{y\})\}}
=1_{\txt\{\Xb_{s,t}(x)\geq z\mbox{ for some }z\in\{y\}\}},
\ee
which by Proposition~\ref{P:Graystar} implies Gray's result
Theorem~\ref{T:Gray}.

\subsection{Cooperative branching}\label{S:coop}

In the present subsection, we consider a more concrete example of an
interacting particle system that is monotone but not additive and to which
Theorem~\ref{T:mondual} applies. We will give a concrete description of the
dual process with generator $H_\ast$ from (\ref{HPoisY}) and make optimal use
of the fact that some of the maps involved are additive.

Denoting elements of $\{0,1\}^n$ as finite words made up of the letters 0 and
1, consider the maps $a,b,c,d,e$ acting on $\{0,1\}^2$, $\{0,1\}^3$,
$\{0,1\}^2$, $\{0,1\}^1$, and $\{0,1\}^2$, respectively, defined as follows.
\be\ba{lll}\label{coomaps}
\mbox{voter move}\quad
&a(10):=00,\ a(01):=11,\quad&a(x):=x\mbox{ otherwise,}\\[5pt] 
\mbox{cooperative branching}\quad
&b(110):=111,\quad&b(x):=x\mbox{ otherwise,}\\[5pt] 
\mbox{coalescing RW}\quad
&c(11):=01,\ c(10):=01,\quad&c(x):=x\mbox{ otherwise,}\\[5pt] 
\mbox{death}\quad
&d(1):=0,\quad&d(x):=x\mbox{ otherwise,}\\[5pt]
\mbox{exclusion}\quad
&e(10):=01,\ e(01):=10,\quad&e(x):=x\mbox{ otherwise.}
\ec
These maps can be lifted to a larger space of the form $\{0,1\}^\La$ with
$\La$ any finite set as follows. For each ordered triple $(i,j,k)\in\La^3$
such that $i,j,k\in\La$ are all different from each other, let $x(ijk)$ denote
the word obtained by writing $x(i),x(j)$, and $x(k)$ after each other, and let
$b_{ijk}:\{0,1\}^\La\to\{0,1\}^\La$ be the map defined by
\be\ba{l}\label{lift}
\dis\big(b_{ijk}x\big)(ijk):=b\big(x(ijk)\big)\quand
\big(b_{ijk}x\big)(m):=x(m)\quad\mbox{for}\quad m\not\in\{i,j,k\}.
\ec
We define $a_{ij},c_{ij},d_i$, and $e_{ij}$ in a similar way, by lifing the map
$a,c,d$, and $e$ to $\{0,1\}^\La$.

The map $a_{ij}$ corresponds to the voter model map ${\rm vot}_{ji}$ (note the
order of the indices) from (\ref{votmap}) and the map $c_{ij}$ coincides with
the coalescing random walk map ${\rm rw}_{ij}$ from (\ref{rwmap}). The death
map $d_i$ corresponds to deaths of particles at $i$ and $e_{ij}$ corresponds
to exclusion process dynamics, which interchanges the states of the sites $i$
and $j$. The map $b_{ijk}$ describes the situation where two particles at
$i,j$ are both needed to produce a third particle at $k$; such maps have been
used to model, e.g., sexual reproduction. In \cite{Nob92,Neu94}, particle
systems whose generators can be represented in the maps $b_{ijk},d_i$, and
$e_{ij}$ are studied, while \cite{SS14} is concerned with systems that involve
the maps $b_{ijk}$ and $c_{ij}$.

We set $S:=\{0,1\}^\La$. Here, we will only consider finite $\La$. The
extension to infinite $\La$ such as $\La=\Z$, considered  in
\cite{Nob92,Neu94,SS14}, will be treated in Section~\ref{S:monpart}.
Choose $S':=S$, $x':=1-x$, so that $S'$ is a dual
of $S$ as in the sense of (\ref{duspa}). With this choice, the function
$\li\,\cdot\,,\,\cdot\,\re$ from (\ref{tetSS}) becomes
\be
\li x,y\re=1_{\txt\{x\leq y'\}}=1_{\txt\{x\wedge y=0\}}\qquad(x,y\in S).
\ee
We start by observing that the maps $a_{ij},c_{ij},d_i$, and $e_{ij}$ are
additive and hence have duals in the sense of Lemma~\ref{L:addual}.

\bl[Additive maps]\label{L:acde}
For each $i,j\in\La$ such that $i\neq j$, the maps $a_{ij},c_{ij},d_i$, and
$e_{ij}$ are additive and their dual maps in the sense of Lemma~\ref{L:addual}
are given by
\be
a'_{ij}=c_{ij},\quad c'_{ij}=a_{ij},\quad d'_i=d_i,\quand e'_{ij}=e_{ij}.
\ee
\el

The maps $b_{ijk}$, on the other hand, are not additive, but only
monotone. Sticking to our choice $S':=S$ and $x':=1-x$, we fall back on the
duality function $\phi$ from (\ref{phidef}), which now reads
\be\label{phidef3}
\phi(x,B)
=1_{\txt\{x\leq y'\mbox{ for some }y\in B\}}
=1_{\txt\{x\wedge y=0\mbox{ for some }y\in B\}}
\ee
$\big(x\in S,\ B\in\Pc(S)\big)$. Our next aim is to determine the dual map
$b^\ast_{ijk}$ defined in (\ref{mast}). We start with the map $b$ from
(\ref{coomaps}). We observe that
\be
\big(b^{-1}\big(\{x\}^\down\big)\big)_{\rm max}
=\left\{\ba{ll}
\{100,010\}\quad&\mbox{if }x=110,\\[5pt]
\{x\}\quad&\mbox{otherwise.}
\ea\right.
\ee
By (\ref{mast}), taking into account that $x':=1-x$, it follows that
\be
b^\ast\big(\{x\}\big)
=\left\{\ba{ll}
\{011,101\}\quad&\mbox{if }x=001,\\[5pt]
\{x\}\quad&\mbox{otherwise.}
\ea\right.
\ee

In order to find a more convenient expression for $b^\ast$, we define maps
$b^{(1)},b^{(2)}:\{0,1\}^3\to\{0,1\}^3$ by
\be\ba{ll}\label{b12}
\dis b^{(1)}(001):=011,\quad&\dis b^{(1)}(x):=x\mbox{ otherwise,}\\[5pt]
\dis b^{(2)}(001):=101,\quad&\dis b^{(2)}(x):=x\mbox{ otherwise.}
\ec
Then, using (\ref{cuprop}), we see that for any $B\sub\{0,1\}^3$,
\be
b^\ast(B)=\{b^{(1)}(x):x\in B\}\cup\{b^{(2)}(x):x\in B\}
=b^{(1)}(B)\cup b^{(2)}(B),
\ee
where $b^{(1)}(B)$ and $b^{(2)}(B)$ denote the images of $B$ under the maps
$b^{(1)}$ and $b^{(2)}$, respectively. Similarly, lifting these maps to the
larger space $S=\{0,1\}^\La$ in the same way as before,
\be\label{bast}
b^\ast_{ijk}(B)=b^{(1)}_{ijk}(B)\cup b^{(2)}_{ijk}(B)
\qquad\big(B\in\Pc(S)\big).
\ee
Lemma~\ref{L:astac} tells us that $a^\ast_{ij}(B)=a'_{ij}(B)$, where
$a'_{ij}(B)$ denotes the image of $B$ under the map $a'_{ij}$,
and similarly for the additive maps $c_{ij},d_i$, and $e_{ij}$.
Using also Lemma~\ref{L:acde}, we see that
\be\label{cast}
a^\ast_{ij}(B)=c_{ij}(B),\quad c^\ast_{ij}(B)=a_{ij}(B),\quad
d^\ast_i(B)=d_i(B),\quand e^\ast_{ij}(B)=e_{ij}(B).
\ee

By Lemma~\ref{L:dumaps}, the maps $a^\ast_{ij},b^\ast_{ij},\ldots$ in
(\ref{bast}) and (\ref{cast}) are dual to $a_{ij},b_{ij},\ldots$ with respect
to the duality function $\phi$ from (\ref{phidef3}), so
Proposition~\ref{P:pathdual} tells us how to construct a pathwise dual for a
Markov process whose generator is represented in these maps.

It is instructive to do this in a concrete example.
Let us consider a process $X$ with a generator of the form
\be\label{Gcoop}
Gf(x):=\sum_{ijk}r_{ijk}\big(f(b_{ijk}(x))-f(x)\big)
+\sum_{ij}s_{ij}\big(f(c_{ij}(x))-f(x)\big),
\ee
where the first sum runs over all ordered triples $(i,j,k)\in\La^3$
such that $i,j,k\in\La$ are all different from each other, the second sum runs
over all $(i,j)\in\La^2$ with $i\neq j$, and $r_{ijk}$ and $s_{ij}$ are
nonnegative rates. The process $X$ is a process of particles performing
cooperative branching and coalescing random walk dynamics such as studied in
\cite{SS14}.

By Proposition~\ref{P:pathdual}, Lemma~\ref{L:dumaps}, (\ref{bast}), and
(\ref{cast}), the process $X$ is pathwise dual with respect to the function
$\phi$ from (\ref{phidef3}) to the $\Pc(S)$-valued process
$Y^\ast=(Y^\ast_t)_{t\geq 0}$ with generator
\be
\label{Gast}
H_\ast f(B):=\sum_{ijk}r_{ijk}
\big(f\big(b^{(1)}_{ijk}(B)\cup b^{(2)}_{ijk}(B)\big)-f(B)\big)
+\sum_{ij}s_{ij}\big(f(a_{ij}(B))-f(B)\big).
\ee

We can think of the set-valued process $(Y^\ast_t)_{t\geq 0}$ as an evolving
collection of voter-model configurations. With rate $s_{ij}$, the voter map
$a_{ij}$ is applied to all configurations $y\in Y^\ast_t$ simultaneously. With
rate $r_{ijk}$, we apply both of the maps $b^{(1)}_{ijk}$ and $b^{(2)}_{ijk}$
to each configuration $y\in Y^\ast_t$, and collect all different outcomes in a
new collection of voter-model configurations.

Since each pathwise dual is also a (normal) dual, we conclude that
if $X$ and $Y$ are processes with generators as in (\ref{Gcoop}) and
(\ref{Gast}) deterministic initial states $X_0$ and $Y_0$, then
\be
\P\big[X_t\wedge y=0\mbox{ for some }y\in Y_0\big]
=\P\big[X_0\wedge y=0\mbox{ for some }y\in Y_t\big]
\qquad(t\geq 0),
\ee
or equivalently
\be
\P\big[X_t\wedge y\neq 0\ \forall y\in Y_0\big]
=\P\big[X_0\wedge y\neq 0\ \forall y\in Y_t\big]
\qquad(t\geq 0).
\ee
Using (\ref{cast}), it is straightforward to extend this duality so
that the generator of $X$ also includes application of the maps $d_i$ and
$e_{ij}$ with certain rates as in \cite{Nob92,Neu94}.

\section{Proofs}\label{S:proofs}

\subsection{Overview}

In this section, we provide the proofs for all facts stated so far without
proof or reference. Proposition~\ref{P:pathdual}, Lemma~\ref{L:genrep}, and
Lemma~\ref{L:invdual} are proved in Section~\ref{S:Markdual}.
Lemma~\ref{L:monadd} is proved in Section~\ref{S:monad} while
Lemma~\ref{L:addual} has already been proved ``on the spot'' in
Section~\ref{S:addual}. Theorem~\ref{T:addual} follows directly from
Proposition~\ref{P:pathdual} and Lemma~\ref{L:addual}.
Section~\ref{S:monproof} contains the proofs of Lemma~\ref{L:dumaps} and the
later Lemma~\ref{L:astac}. The preceding Theorem~\ref{T:mondual} follows
directly from Proposition~\ref{P:pathdual} and Lemma~\ref{L:dumaps}, and, as
already pointed out in the text, Lemma~\ref{L:altmaps} follows by applying
Lemma~\ref{L:dumaps} to the reversed order. Proposition~\ref{P:monsuf} has
been cited in the text from \cite{KKO77,FM01}. Lemmas~\ref{L:madd},
\ref{L:percrep}, and \ref{L:nondis} are proved in
Section~\ref{S:percproof}. As explained in the text preceding it,
Theorem~\ref{T:Siegmund} follows from Theorem~\ref{T:addual}. It has also
already been explained that Proposition~\ref{P:Siegmund} follows from
Theorem~\ref{T:addual}. Lemma~\ref{L:Krodu} is proved in
Section~\ref{S:ducalc}, which also contains the proof of the later
Lemma~\ref{L:acde}, while the preceding Lemma~\ref{L:spinmon} and
Proposition~\ref{P:Graystar} are proved in Section~\ref{S:Grayproof}.
Theorem~\ref{T:Gray}, finally, is cited from \cite{Gra86} and as explained in
the text also follows from Proposition~\ref{P:pathdual},
Lemma~\ref{L:altmaps}, and Proposition~\ref{P:Graystar}.

\subsection{Markov process duality}\label{S:Markdual}

In this section we prove Lemma~\ref{L:Markrep} and Proposition~\ref{P:pathdual}
about the Poisson construction of Markov processes with finite state spaces
and pathwise duality, as well as Lemma~\ref{L:genrep} about representability of
Markov generators and semigroups and the simple observation
Lemma~\ref{L:invdual} about the inverse image map.

\begin{Proof}[of Proposition~\ref{P:pathdual}]
If
\be
\De_{t-,u-}=\{(m_1,t_1),\ldots,(m_n,t_n)\},
\qquad t_1<\cdots<t_n,
\ee
then
\be
\Xb_{t-,u-}:=m_n\circ\cdots\circ m_1
\quand
\Yb_{-u,-t}:=\hat m_1\circ\cdots\circ\hat m_n,
\ee
so repeated application of (\ref{mdual}) implies that 
$\Xb_{t-,u-}$ is dual to $\Yb_{-u,-t}$. Using the semigroup property,
it follows that the expression
\be
\psi\big(\Xb_{s,t-}(x),\Yb_{-u,-t}(y)\big)
=\psi\big(\Xb_{t-,u-}\circ\Xb_{s-,t-}(x),y\big)
=\psi\big(\Xb_{s,u-}(x),y\big)
\ee
is constant as a function of $t\in[s,u]$. For any $t_0<\cdots<t_n$, the pairs
\[
(\Xb_{t_0,t_1},\Yb_{-t_1,-t_0}),\ldots,(\Xb_{t_{n-1},t_n},\Yb_{-t_n,-t_{n-1}})
\]
are functions of the restrictions of the Poisson point process $\De$
to disjoint sets $\Gi\times(t_0,t_1]$ etc., and therefore independent,
completing the proof that the stochastic flows $(\Xb_{s,t})_{s\leq t}$ and
$(\Yb_{s,t})_{s\leq t}$ are dual in the sense defined in
Section~\ref{S:pathdual}. In particular, this proves that the Markov
processes $X$ and $Y$ are pathwise dual as defined in
Section~\ref{S:pathdual}.
\end{Proof}

\begin{Proof}[of Lemma~\ref{L:genrep}]
If $G$ can be represented in $\Gi$, then we can construct a stochastic flow
$(\Xb_{s,t})_{s\leq t}$ based on such a random mapping representation as in
(\ref{Phidef}). By the fact that $\Gi$ is closed under composition and
contains the identity map, it follows that $\Xb_{s,t}\in\Gi$ for each $s\leq
t$, so (\ref{PXb}) proves that $P_t$ can be represented in $\Gi$ for all
$t\geq 0$.

Assume, conversely, that $P_t$ can be represented in $\Gi$ for all
$t\geq 0$. Then, for each $t\geq 0$, we can find a probability distribution
$\pi_t$ on $\Gi$ such that
\be
P_t(x,y)=\sum_{m\in\Gi}\pi_t(m)1_{\txt\{m(x)=y\}}\qquad(x,y\in S).
\ee
Since $G$ is the generator of $(P_t)_{t\geq 0}$,
\be\label{piG}
\lim_{t\down0}t^{-1}\sum_{m\in\Gi}\pi_t(m)1_{\txt\{m(x)=y\}}=G(x,y)
\qquad(x,y\in S,\ x\neq y).
\ee
Let $m\in\Gi$ satisfy $m\neq 1$, i.e., $m$ is different from the identity
map. Then we can find $x\neq y$ such that $m(x)\neq y$. Now (\ref{piG}) shows
that $t^{-1}\pi_t(m)\leq G(x,y)$, so using compactness we can find a sequence
of times $t_n\down 0$ such that the limit
\be
r_m:=\lim_{n\to\infty}t_n^{-1}\pi_{t_n}(m)
\ee
exists for every $m\in\Gi':=\{m\in\Gi:m\neq 1\}$. Using (\ref{piG}), it
follows that
\be
G(x,y)=\sum_{m\in\Gi'}r_m1_{\txt\{m(x)=y\}}\qquad(x,y\in S,\ x\neq y).
\ee
Using also that $G$ is a Markov generator, this is easily seen to imply
(\ref{GPois}), i.e., we have found a random mapping representation of $G$ in
terms of maps in $\Gi$.
\end{Proof}

\begin{Proof}[of Lemma~\ref{L:invdual}]
This is immediate from the observation that
\[
\psi\big(m(x),A\big)=1_{\txt\{m(x)\in A\}}
=1_{\txt\{x\in m^{-1}(A)\}}=\psi\big(x,m^{-1}(A)\big).
\]
\end{Proof}

\subsection{Monotone and additive maps}\label{S:monad}

In this section we prove Lemma~\ref{L:monadd}, which characterizes montone and
additive maps in terms of their inverse images.\med

\begin{Proof}[of Lemma~\ref{L:monadd}]
{\rm (i)} If $S$ and $T$ are partially ordered sets, $m:S\to T$ is monotone,
and $A\sub T$ is decreasing, then for each $x,y\in S$ with $x\leq y$ and $y\in
m^{-1}(A)$, we have $m(x)\leq m(y)\in A$ by the fact that $m$ is monotone and
hence $m(x)\in A$ or equivalently $x\in m^{-1}(A)$ by the fact that $A$ is
decreasing. This shows that
\be
\label{monequiv}
 m^{-1}(A)\in\Pc_{\rm dec}(S) \textrm{ for all }A\in\Pc_{\rm dec}(T).
\ee
Conversely, if (\ref{monequiv}) holds, and $x,y\in S$ satisfy $x\leq y$, then
$m^{-1}(\{y\}^\down)$ is a decreasing set containing $y$, so $x\in
m^{-1}(\{y\}^\down)$ or equivalently $m(x)\leq m(y)$, proving that $m$ is
monotone.

{\rm (ii)} We use the fact that if $S$ is a join-semilattice, then a set
$A\sub S$ is an ideal if and only if $A$ is nonempty, decreasing, and $x,y\in
A$ imply $x\vee y$ in $A$ (see Lemma~\ref{L:latfil} in the appendix).  Now let
$S$ and $T$ be join-semilattices that are bounded from below. Assume that
$m:S\to T$ is additive and that $A\in\Pc_{!\rm dec}(T)$ is an ideal.  Since
additive functions are monotone, by what we have just proved, $m^{-1}(A)$ is
decreasing. Since $A$ is nonempty and decreasing, it contains $0$, and since
$m(0)=0$ we have $0\in m^{-1}(A)$, which shows that $m^{-1}(A)$ is
nonempty. Finally, if $x,y\in m^{-1}(A)$, then by additivity and the fact that
$A$ is an ideal, $m(x\vee y)=m(x)\vee m(y)\in A$, so $x\vee y\in m^{-1}(A)$,
which completes the proof that $m^{-1}(A)$ is an ideal.

Assume, conversely, that $m:S\to T$ has the property that $m^{-1}(A)$ is an
ideal whenever $A\sub T$ is an ideal. Then, since $m^{-1}(\{0\})$ is nonempty
and decreasing, we have $0\in m^{-1}(\{0\})$ and hence $m(0)=0$. For each
$y\in S$, the set $m^{-1}(\{m(y)\}^\down)$ is a decreasing set containing $y$,
so $x\leq y$ implies $x\in m^{-1}(\{m(y)\}^\down)$ or equivalently $m(x)\leq
m(y)$, showing that $m$ is monotone. Since $x\leq x\vee y$ and $y\leq x\vee
y$, this implies that $m(x)\vee m(y)\leq m(x\vee y)$ for each $x,y\in S$. To
get the opposite inequality, we observe that $m^{-1}(\{m(x)\vee m(y)\}^\down)$
is an ideal containing $x$ and $y$, so using the properties of ideals also
$x\vee y\in m^{-1}(\{m(x)\vee m(y)\}^\down)$ which says that $m(x\vee y)\leq
m(x)\vee m(y)$.
\end{Proof}

\subsection{Monotone systems duality}\label{S:monproof}

In this section we prove Lemmas~\ref{L:dumaps} and \ref{L:astac} that form the
basis of Section~\ref{S:mondual} about monotone systems duality.

To prepare for the proof of Lemma~\ref{L:dumaps}, we make the following
general observation on sets of maximal and minimal elements. If $A \sub S$
is finite, then there exists to any $x \in A$ a maximal element $y \in A$ with
$y \geq x$ (see Lemma~\ref{L:minmax} in the appendix). Therefore, we have
$(A_{\rm max})^\down\supset A$ and similarly $(A_{\rm min})^\up\supset A$ for
each $A\in\Pc(S)$, so in particular
\be\label{maxup}
(A_{\rm max})^\down=A\quad\big(A\in\Pc_{\rm dec}(S)\big)
\quand
(A_{\rm min})^\up=A\quad\big(A\in\Pc_{\rm inc}(S)\big).
\ee
(In particular, this also holds if $A=\emptyset$ since $\emptyset_{\rm
  min}=\emptyset=\emptyset_{\rm max}$ and
$\emptyset^\up=\emptyset=\emptyset^\down$.) In view of (\ref{maxup}), we can
``encode'' a set $A\in\Pc_{\rm dec}(S)$ by describing its set of maximal
elements, or more generally by any subset $B\sub S$ such that $A=B^\down$.
\med

\begin{Proof}[of Lemma ~\ref{L:dumaps}]
Through the bijection $x\mapsto x'$, the monotone map $m:S\to S$ naturally
gives rise to a monotone map $n:S'\to S'$ defined by
\be\label{nacm}
n(x):=m(x')'\qquad(x\in S').
\ee
In terms of the map $n$, the definitions (\ref{mast}) take the simpler form
\be\label{mastn}
m^\dgg(B):=n^{-1}(B^\up)_{\rm min}
\quand
m^\ast(B):=\bigcup_{x\in B}n^{-1}(\{x\}^\up)_{\rm min}
\ee
$(B\in\Pc(S'))$. Instead of showing that $m^\dgg$ and $m^\ast$ are dual to $m$
with respect to the duality function in (\ref{phidef}), we may equivalently
show that $m^\dgg$ and $m^\ast$ are dual to $n$ with respect to the duality
function
\be
\ti\phi(x,B):=1_{\txt\{x'\in{B'}^\down\}}=1_{\txt\{x\in B^\up\}}
\qquad\big(x\in S',\ B\in\Pc(S')\big).
\ee
Letting $\both=\dgg$ or $\ast$, this means that we must show that
\be\label{ovmdu}
n(x)\in B^\up\quad\mbox{if and only if}\quad
x\in m^\both(B)^\up
\qquad\big(x\in S',\ B\in\Pc(S')\big).
\ee
Since the event in the left-hand side of this equation is $\{x\in
n^{-1}(B^\up)\}$, this is equivalent to
\be\label{phidu}
m^\both(B)^\up=n^{-1}(B^\up)\qquad\big(B\in\Pc(S')\big).
\ee
Since $n$ is monotone, Lemma~\ref{L:monadd} tells us that
$n^{-1}(A)\in\Pc_{\rm inc}(S')$ for each $A\in\Pc_{\rm inc}(S')$.
Therefore, by (\ref{maxup}) applied to the increasing set $A=n^{-1}(B^\up)$,
we see that
\be
m^\dgg(B)^\up=\big(n^{-1}(B^\up)_{\rm min}\big)^\up=m^{-1}(B^\up)
\qquad\big(B\in\Pc(S')\big),
\ee
proving that $m^\dgg$ satisfies (\ref{phidu}). Similarly, for the map
$m^\ast$, applying (\ref{maxup}) to the increasing sets $A=n^{-1}(\{x\}^\up)$
with $x\in B$, we have for $B\in\Pc(S')$
\be
m^\ast(B)^\up
=\big(\bigcup_{x\in B}n^{-1}(\{x\}^\up)_{\rm min}\big)^\up
=\bigcup_{x\in B}\big(n^{-1}(\{x\}^\up)_{\rm min}\big)^\up
=\bigcup_{x\in B} n^{-1}(\{x\}^\up)
=n^{-1}(B^\up),
\ee
which shows that also $m^\ast$ satisfies (\ref{phidu}).

It remains to prove the properties (\ref{dagprop}) and (\ref{cuprop}).
It is clear from (\ref{mastn}) that
\be\label{dagcup}
m^\dgg(B)=m^\dgg(B)_{\rm min}
\quand
m^\ast(B\cup C)=m^\ast(B)\cup m^\ast(C)
\qquad\big(B,C\in\Pc(S')\big).
\ee
Formula (\ref{phidu}) tells us that
\be
m^\dgg(B)^\up=n^{-1}(B^\up)=m^\ast(B)^\up,
\ee
Taking minima on both sides, using the facts that $m^\dgg(B)=m^\dgg(B)_{\rm
  min}$ and $(A^\up)_{\rm min}=A_{\rm min}$ for any $A\sub S$, we see that
$m^\dgg(B)=(m^\ast(B))_{\rm min}$ for $B\in\Pc(S')$.
\end{Proof}

\begin{Proof}[of Lemma~\ref{L:astac}]
By Lemma~\ref{L:addual}, we have
\be\ba{l}
\dis m^{-1}(\{y'\}^\down)
=\{x\in S:m(x)\leq y'\}
=\{x\in S:\li m(x),y\re=1\}\\[5pt]
\dis\quad=\{x\in S:\li x,m'(y)\re=1\}
=\{m'(y)'\}^\down\qquad(y\in S').
\ec
It follows that $(m^{-1}(\{y'\}^\down))_{\rm max}=\{m'(y)'\}$, so by
(\ref{mast}), we see that $m^\ast(B)=\{m'(y):y\in B\}$ $(B\in\Pc(S))$.
The fact that $m^\dgg(B)=m^\ast(B)_{\rm min}$ has already been proved in
(\ref{dagprop}).
\end{Proof}

\subsection{Percolation representations}\label{S:percproof}

In this section we prove Lemmas~\ref{L:madd}, \ref{L:percrep}, and
\ref{L:nondis}, that show that each additive system taking values in a lattice
has a percolation representation, and if the lattice is moreover distributive,
then the process and its dual can be represented in the same percolation
substructure.\med

\begin{Proof}[of Lemma~\ref{L:madd}]
We start by showing that if $M\sub \La\times\La$ satisfies (\ref{Mprop}),
then (\ref{minM}) defines an additive map $m:\Pc_{\rm dec}(\La)\to\Pc_{\rm
  dec}(\La)$. To see that $m$ maps $\Pc_{\rm dec}(\La)$ into itself, we note
that $\ti\jj\leq j$ and $j \in m(x)$ implies $(i, j)\in M$ for some $ i \in x$
which by (\ref{Mprop})~(ii) implies that $(i,\ti\jj)\in M$ and hence
$\ti\jj\in m(x)$. We see from (\ref{minM}) that
\be
m(x)=\bigcup_{i \in x}m\big(\{i\}\big)\qquad\big(x\in\Pc_{\rm dec}(\La)\big),
\ee
which is easily seen to imply that $m$ is additive.

We next show that if $m$ is defined via $M$ as in (\ref{minM}), then we can
recover $M$ from $m$ since
\be\label{Mndefproof}
M=\big\{(i,j)\in \La \times \La :j \in m\big(\{i\}^\down\big)\big\}.
\ee
To see this, note that $(i,j)\in M$ implies by (\ref{minM}) $j\in
m\big(\{i\}^\down\big)$ since $i \in \{i\}^\down$.  This proves the inclusion
$\sub$ in (\ref{Mndefproof}). Conversely, if $j\in m\big(\{i\}^\down\big)$,
then by (\ref{minM}) there must exist some $\ti\ii\leq i$ such that
$(\ti\ii,j)\in M$, which by (\ref{Mprop})~(i) implies $(i,j)\in M$ and shows
$\supset$ in (\ref{Mndefproof}).

To see that every additive map from $\Pc_{\rm dec}(\La)$ into itself is of the
form (\ref{minM}), let $m$ be such a map. We will show that $m$ is of the form
(\ref{minM}) where $M$ is given by (\ref{Mndefproof}). Indeed, if $M$ is given
by (\ref{Mndefproof}), then $M$ satisfies (\ref{Mprop})~(i) by the
monotonicity of $m$, namely if $(i,j)\in M$ and $i\leq\ti\ii$ then $j\in
m\big(\{i\}^\down\big)\sub m\big(\{\ti\ii\}^\down\big)$, which implies
$(\ti\ii,j)\in M$. Property (\ref{Mprop})~(ii) holds due to the fact that
$m\big(\{i\}^\down\big)\in\Pc_{\rm dec}(\La)$ so $j\in m\big(\{i\}^\down\big)$ implies $\ti\jj
\in m\big(\{i\}^\down\big)$ for any $\ti\jj\leq j.$

To see that $m$ is given by (\ref{minM}), first define $n$ by (\ref{minM}). We
claim that for any $i\in\La$,
\be
n\big(\{i\}^\down\big)
=\big\{j:( \ti\ii,j)\in M\mbox{ for some } \ti\ii \leq i\big\}
=\big\{j:(i,j)\in M\big\}=m\big(\{i\}^\down\big).
\ee
Indeed, in the second equality, we have $\sub$ by (\ref{Mprop})~(i) and
$\supset$ by choosing $\ti\ii=i$, while the third equality is immediate from
(\ref{Mndefproof}). It follows that $m=n$ on sets of the form
$\{x\}^\down$. Since both $m$ and $n$ are additive, it follows that they agree
on $\Pc_{\rm dec}(\La)$.

It remains to prove (\ref{Madj}). Recall that $S:=\Pc_{\rm dec}(\La)$ and that
$S':=\Pc_{\rm dec}(\La')=\Pc_{\rm inc}(\La)$, which together with the map $x\mapsto
x':=x^{\rm c}$ is a dual to $S$. Then the duality function takes the form
$\li x,y\re=1_{\txt\{x\cap y=\emptyset\}}$ and two additive maps $m:S\to S$ and 
$m':S'\to S'$ are dual if and only if
\be\label{mmac}
1_{\txt\{m(x)\cap y=\emptyset\}}=1_{\txt\{x\cap m'(y)=\emptyset\}}
\qquad\big(x\in\Pc_{\rm dec}(\La),\ y\in\Pc_{\rm inc}(\La)\big).
\ee
We observe that if $M$ satisfies (\ref{Mprop}) and $M'$ is given by
(\ref{Madj}), then $M'$ satisfies (\ref{Mprop}) with respect to the reversed
order, i.e., with respect to the order on $\La'$. In view of this, it suffices
to show that if $m:S\to S$ is defined in terms of $M$ as in (\ref{minM}) and
similarly $m':S'\to S'$ is defined in terms of $M'$, then $m$ and $m'$ are
dual in the sense of (\ref{mmac}). Indeed,
\be\ba{r@{\quad}c@{\quad}l@{\quad}l@{\quad}l}
m(x)\cap y=\emptyset
&\desd&
\{j\in\La:(i,j)\in M\mbox{ for some }i\in x\}\cap y=\emptyset\\[5pt]
&\desd&
\{(i,j)\in M:i\in x,\ j\in y\}=\emptyset\\[5pt]
&\desd&
x\cap\{i\in\La:(j,i)\in M'\mbox{ for some }j\in y\}=\emptyset
&\desd&
x\cap m'(y)=\emptyset.
\ec
\end{Proof}

\begin{Proof}[of Lemma~\ref{L:percrep}]
The left- and right-hand side of (\ref{percrep}) both equal $x$ when $s=u$ and
for fixed $s\in\R$, both sides change only at times $u$ when $(m,u)\in\De$ for
some $m\in\Gi$. If just before time $u$, both sides are equal to some
$y\in\Pc_{\rm dec}(\La)$, then at time $u$, using the definition of open
paths, we see that the left- and right-hand side of (\ref{percrep}) equal
\be
m(y)\quand\{j\in\La:(i,j)\in M\mbox{ for some }i\in y\},
\ee
respectively. In view of this, (\ref{percrep}) follows by induction from
(\ref{minM}). The proof of (\ref{dualrep}) is the same: we first observe that
equality holds at $u=s$ and then check that the equation remains true when we
increase $u$ while keeping $s$ fixed, where
\be
m'(y)=\{j\in\La:(i,j)\in M'\mbox{ for some }i\in y\},
\ee
with $M'$ as in (\ref{Madj}).

Let us write $(i,s)\leadsto(j,u-)$ if there is an open path ``from $(i,s)$ to
$(j,u-)$'', where in the definition of an open path, we replace the time
interval $(s,u]$ by $(s,u)$. In the same way, we can give a meaning to
$(i,s-)\leadsto(j,u)$ and $(i,s-)\leadsto(j,u-)$, where arrows and blocking
symbols at the end of the time interval have an effect, or not, depending on
whether we want our definition to be left- or right-continuous in $s$ or $u$.
Then formulas (\ref{percrep}) and (\ref{dualrep}) have obvious analogues where
right-continuity in a variable is replaced by left-continuity, or vice versa.
Now (\ref{duint}) follows from the observation that
\bc
\dis\Xb_{s,t-}(x)\cap\Yb_{-u,-t}(y)=\emptyset
&\quad\desd\quad&\dis
\not\exists i,j,k\mbox{ such that }i\in x,\ j\in y,
\ (i,s)\leadsto(k,t-)\leadsto(j,u-)\\[5pt]
&\quad\desd\quad&\dis
(i,s)\not\leadsto(j,u)\ \forall\ i\in x,\ j\in y\quad{\rm a.s.},
\ec
where in the last step we have used a.s.\ continuity at the deterministic
time~$u$.
\end{Proof}

\begin{Proof}[of Lemma~\ref{L:nondis}]
Let $S$ be a finite lattice. Then we claim that the
map\footnote{Alternatively, one may consider the map
  $x\mapsto\big(\{x\}^\up\cap\La\big)^{\rm c}$ where $\La$ is the set of
  meet-irreducible elements of $S$. As a result of Birkhoff's representation
  theorem, one can prove that this map is onto, and as a result sets up an
  isomorphism between the lattices $S$ and $\Pc_{\rm dec}(\La)$, if and only
  if $S$ is distributive.}
\be\label{joinisom}
x\mapsto\big(\{x\}^\up\big)^{\rm c}
\ee
is a $(0,\vee)$-homomorphism from $S$ into the lattice of sets $\Pc_{\rm
  dec}(S)$. Indeed,
\be
\big(\{0\}^\up\big)^{\rm c}=\emptyset
\quand
\big(\{x\vee y\}^\up\big)^{\rm c}
=\big(\{x\}^\up\cap\{y\}^\up\big)^{\rm c}
=\big(\{x\}^\up\big)^{\rm c}\cup\big(\{y\}^\up\big)^{\rm c}.
\ee
Since $x\neq y$ implies $\{x\}^\up\neq\{y\}^\up$, the map in (\ref{joinisom})
is one-to-one and as a result a $(0,\vee)$-isomorphism to its image, proving
that each finite lattice is $(0,\vee)$-isomorphic to a join-semilattice of
sets.

Now let $\La$ be a finite set and let $T\sub\Pc(\La)$ be a join-semilattice of
sets, i.e., $\emptyset\in T$ and $T$ is closed under unions. We claim that
each $m:T\to T$ can be extended to an additive map $\ov
m:\Pc(\La)\to\Pc(\La)$. We will actually prove a somewhat stronger statement.
Assume moreover that $\La$ is partially ordered and that $T\sub\Pc_{\rm
  dec}(\La)$. Then we will show that $m:T\to T$ can be extended to an additive
map $\ov m:\Pc_{\rm dec}(\La)\to\Pc_{\rm dec}(\La)$. In particular, equipping
$\La$ with the trivial order, this includes the statement in
Lemma~\ref{L:nondis} as a special case. The following argument was suggested
to us by L\'aszl\'o Csirmaz.

Assume that $x\in\Pc_{\rm dec}(\La)$ and $x\not\in T$. Then $\ov T:=\{y,x\cup
y:y\in T\}\sub\Pc_{\rm dec}(\La)$ is $(\emptyset,\cup)$-closed and contains $x$.
By the finiteness of $S$, the claim will follow by induction if we
can prove that $m$ can be extended to an additive map 
$\ov m:\ov T \to \Pc_{\rm dec}(\La)$.

For $y\in T$, we define
\bc\label{mAX}
\dis\ov m(y)&:=&m(y),\\[5pt]
\dis\ov m(x\cup y)
&:=&\dis m(y)\cup\bigcap\{m(z):x\sub z\in T\}\\[5pt]
&=&\dis\bigcap\{m(y)\cup m(z):x\sub z\in T\}\\[5pt]
&=&\dis\bigcap\{m(y\cup z):x\sub z\in T\}.
\ec
Note that since $\Pc_{\rm dec}(\La)$ is closed under finite intersections, the last
line in (\ref{mAX}) defines an element of $\Pc_{\rm dec}(\La)$.
We need to show that (\ref{mAX}) is a good definition, i.e.,
\begin{enumerate}
\item If $y\in T$ and $x\cup y\in T$, then
 $\dis\bigcap\{m(y\cup z):x\sub z\in T\}=m(x\cup y)$.
\item If $y,y'\in T$ and $x\cup y=x\cup y'$, then\\[5pt]
 $\dis\bigcap\{m(y\cup z):x\sub z\in T\}=\bigcap\{m(y'\cup z):x\sub z\in T\}$.
\end{enumerate}
Indeed, in (i), the inclusion $\supset$ follows from the monotonicity of $m$,
while the inclusion $\sub$ follows by setting $z=x\cup y$. Property~(ii)
follows from the observation that
\be
x\cup y=x\cup y'\quand x\sub z\in T
\quad\mbox{imply}\quad y\cup z=y'\cup z.
\ee
To show that $\ov m$ is additive, we must prove that
\begin{enumerate}
\item If $y,y'\in T$, then $\ov m(y\cup y')=\ov m(y)\cup \ov m(y')$.
\item If $y,y'\in T$, then $\ov m\big((x \cup y)\cup y'\big)
      =\ov m(x \cup y)\cup\ov m(y')$.
\item If $y,y'\in T$, then $\ov m\big((x \cup y)\cup(y'\cup x)\big)
      =\ov m(x \cup y)\cup\ov m(y'\cup x)$.
\end{enumerate}
Property~(i) follows from the additivity of $m$. Properties~(ii) and (iii) say
that for $y,y'\in T$,
\be
\ov m(y\cup y'\cup x)=\ov m(x \cup y)\cup\ov m(y')
=\ov m(x \cup y)\cup\ov m(x\cup y'),
\ee
which is easily seen to follow from our definition $\ov m(x \cup y)
:=m(y)\cup\bigcap\{m(z):x\sub z\in T\}$  as well as from the additivity of $m$.
\end{Proof}

\subsection{Examples of dual maps}\label{S:ducalc}

In this section we prove Lemmas~\ref{L:Krodu} and \ref{L:acde} by verifying
that the maps there are really dual to each other as claimed.\med

\begin{Proof}[of Lemma~\ref{L:Krodu}]
The easiest way to see this is from the graphical representation (see
Figure~\ref{fig:Krone}). First, one represents the maps $a_i,b_{ij},c_i,d_i$,
and $e_i$ by arrows and blocking symbols as described at the end of
Section~\ref{S:Krone}, and checks that this really corresponds to the maps in
(\ref{Kromaps}) with the interpretation of the two-stage contact process as a
set-valued process in the forward time direction (\ref{forwrep}). The dual
maps $a'_i,b'_{ij},c'_i,d'_i$, and $e'_i$ can then be found according to the
recipe ``reverse the arrows, keep the blocking symbols'' (as is clear from
Figure~\ref{fig:Krone} and which also corresponds to (\ref{Madj})).  Using the
interpretation of the dual set-valued process (\ref{backrep}), this yields
(\ref{Krodu}).
\end{Proof}

\begin{Proof}[of Lemma~\ref{L:acde}]
As in the proof of Lemma~\ref{L:Krodu}, we use the graphical representation.
The maps $a_{ij},c_{ij},d_i$, and $e_{ij}$ can be represented in terms of
arrows and blocking symbols as:

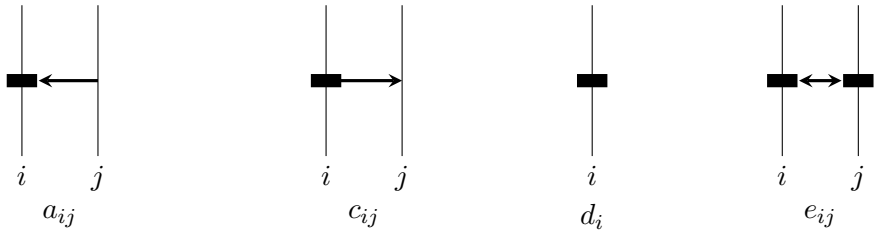
\begin{center}
\begin{tikzpicture}[>=stealth,scale=1]
\begin{scope}
\foreach \x in {1,2}
 \draw (\x,0)--(\x,2);
\draw (1,0) node[below] {$i$};
\draw (2,0) node[below] {$j$};
\draw[very thick,->] (2,1) -- ++ (-0.78,0);
\draw[line width=5pt] (1,1)++(-0.2,0) -- ++ (0.4,0);
\draw (1.5,-0.8) node {$a_{ij}$};
\end{scope}
\begin{scope}[xshift=4cm]
\foreach \x in {1,2}
 \draw (\x,0)--(\x,2);
\draw (1,0) node[below] {$i$};
\draw (2,0) node[below] {$j$};
\draw[very thick,->] (1,1) -- ++ (1,0);
\draw[line width=5pt] (1,1)++(-0.2,0) -- ++ (0.4,0);
\draw (1.5,-0.8) node {$c_{ij}$};
\end{scope}
\begin{scope}[xshift=7.5cm]
\draw (1,0)--(1,2);
\draw (1,0) node[below] {$i$};
\draw[line width=5pt] (1,1)++(-0.2,0) -- ++ (0.4,0);
\draw (1,-0.8) node {$d_i$};
\end{scope}
\begin{scope}[xshift=10cm]
\foreach \x in {1,2}
 \draw (\x,0)--(\x,2);
\draw (1,0) node[below] {$i$};
\draw (2,0) node[below] {$j$};
\draw[very thick,<->] (1.22,1) -- ++ (0.56,0);
\draw[line width=5pt] (1,1)++(-0.2,0) -- ++ (0.4,0);
\draw[line width=5pt] (2,1)++(-0.2,0) -- ++ (0.4,0);
\draw (1.5,-0.8) node {$e_{ij}$};
\end{scope}
\end{tikzpicture}
\end{center}

\noi
According to the recipe ``reverse the arrows, keep the blocking symbols'',
the dual maps are given by:

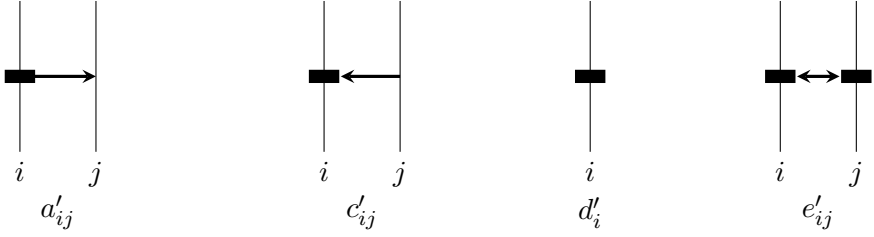
\begin{center}
\begin{tikzpicture}[>=stealth,scale=1]
\begin{scope}
\foreach \x in {1,2}
 \draw (\x,0)--(\x,2);
\draw (1,0) node[below] {$i$};
\draw (2,0) node[below] {$j$};
\draw[very thick,->] (1,1) -- ++ (1,0);
\draw[line width=5pt] (1,1)++(-0.2,0) -- ++ (0.4,0);
\draw (1.5,-0.8) node {$a'_{ij}$};
\end{scope}
\begin{scope}[xshift=4cm]
\foreach \x in {1,2}
 \draw (\x,0)--(\x,2);
\draw (1,0) node[below] {$i$};
\draw (2,0) node[below] {$j$};
\draw[very thick,->] (2,1) -- ++ (-0.78,0);
\draw[line width=5pt] (1,1)++(-0.2,0) -- ++ (0.4,0);
\draw (1.5,-0.8) node {$c'_{ij}$};
\end{scope}
\begin{scope}[xshift=7.5cm]
\draw (1,0)--(1,2);
\draw (1,0) node[below] {$i$};
\draw[line width=5pt] (1,1)++(-0.2,0) -- ++ (0.4,0);
\draw (1,-0.8) node {$d'_i$};
\end{scope}
\begin{scope}[xshift=10cm]
\foreach \x in {1,2}
 \draw (\x,0)--(\x,2);
\draw (1,0) node[below] {$i$};
\draw (2,0) node[below] {$j$};
\draw[very thick,<->] (1.22,1) -- ++ (0.56,0);
\draw[line width=5pt] (1,1)++(-0.2,0) -- ++ (0.4,0);
\draw[line width=5pt] (2,1)++(-0.2,0) -- ++ (0.4,0);
\draw (1.5,-0.8) node {$e'_{ij}$};
\end{scope}
\end{tikzpicture}
\end{center}

\noi
Comparing the pictures, we see that
\[
a'_{ij}=c_{ij},\quad c'_{ij}=a_{ij},\quad d'_i=d_i,\quand e'_{ij}=e_{ij}.
\]
\end{Proof}

\subsection{Gray's duality}\label{S:Grayproof}

In this section, we prove Proposition~\ref{P:Graystar}, which shows that
Gray's dual for general monotone spin systems is a special case of the 
dual $Y^\dott$ from Section~\ref{S:mondual}. For completeness, we also prove
Lemma~\ref{L:spinmon} which says that attractive spin systems are
monotonically representable. We start with a simple preparatory lemma.

\bl[Monotone real functions]\label{L:Fmon}
Let $S$ be a finite partially ordered set. Then $f\in\Fi_{\rm mon}(S,\R)$ if and
only if $f$ can be written as
\be\label{fsum}
f(x)=\sum_{A\in\Pc_{\rm inc}(S)}r_A1_{\{x\in A\}}\qquad(x\in S)
\ee
for real constants $(r_A)_{A\in\Pc_{\rm inc}(S)}$ satisfying $r_A\geq 0$ for
all $A\neq\emptyset,S$.
\el
\begin{Proof}
If $f\in\Fi(S,\R)$ is monotone, then\footnote{Let $S$ and $T$ be partially
  ordered sets. If a map $m:S\to T$ is monotone, then it is
  also monotone with respect to the reversed orders on $S$ and $T$. In view of
  this and Lemma~\ref{L:monadd}~(i), a map $m:S\to T$ is monotone if and only
  if $m^{-1}(A)\in\Pc_{\rm inc}(S)$ for all $A\in\Pc_{\rm inc}(T)$.}
$f^{-1}([r,\infty))\in\Pc_{\rm inc}(S)$ for all $r\in\R$. Let
$f(S):=\{f(x):x\in S\}$ be the image of $S$ under $f$. Since $S$ is finite,
so is $f(S)$ and we may write $f(S)=\{r_1,\ldots,r_n\}$ with $r_1<\cdots<r_n$.
The sets $A_k:=\{x\in S:f(x)\geq r_k\}$ are increasing and
\be
f(x)=r_11_{\{x\in S\}}+\sum_{k=2}^n(r_k-r_{k-1})1_{\{x\in A_k\}}.
\ee
This proves that every $f\in\Fi_{\rm mon}(S,\R)$ can be written in the form
(\ref{fsum}). Conversely, if $f$ is of the form (\ref{fsum}), then it is a sum
of monotone functions, so $f\in\Fi_{\rm mon}(S,\R)$.
\end{Proof}

\begin{Proof}[of Lemma~\ref{L:spinmon}]
It suffices to prove that the first term in (\ref{GGray}) is monotonically
representable. The same arguments applied to the second term and $S$ equipped
with the reversed order then prove the general statement.

Since each $\bet_i$ is monotone and nonnegative, by Lemma~\ref{L:Fmon}, for
each $i\in\La$ we can find some set $\Ai_i$ whose elements are increasing,
nonempty subsets of $S$, as well as nonnegative constants
$(r_{i,A})_{A\in\Ai_i}$, such that
\be
\bet_i(x)=\sum_{A\in\Ai_i}r_{i,A}1_{\txt\{x\in A\}}.
\ee
For each $i\in\La$ and $A\in\Ai_i$, define a map $m_{i,A}$ by
\be
m_{i,A}(x):=\left\{\ba{ll}
x\vee\eps_i\quad\mbox{if }x\in A,\\
x\quad\mbox{otherwise.}\ea\right.
\ee
Then the first term in (\ref{GGray}) can be written as
\be
\sum_{i\in\La}\sum_{A\in\Ai_i}r_{i,A}\big(f(m_{i,A}(x))-f(x)\big),
\ee
which is the desired representation in monotone maps.
\end{Proof}

\begin{Proof}[of Proposition~\ref{P:Graystar}]
For any $s\leq u$ and $B\in\Pc(S)$, let us write
\be
\zeta_{s,u}(B):=\bigcup_{y\in B}\zeta_{s,u}(y).
\ee
We will show that $\Yb^\dott_{-u,-s}(B)=\zeta_{s,u}(B)$. In particular,
setting $B=\{y\}$ then yields the statement in Proposition~\ref{P:Graystar}.
Thus, in light of the definition of Gray's dual in (\ref{zetadef}), we claim
that for any $s\leq u$ and $B\in\Pc(S)$, a.s.
\be\label{Zaim}
\Yb^\dott_{-u,-s}(B)=\big\{z\in S:
\mbox{there exists an $[s,u]$-path $\pi$ from $z$ to some $y\in B$}\big\}.
\ee
As in (\ref{Phidef}), we order the elements of $\De_{s,u}$ according to the
time they occur:
\be
\De_{s,u}:=\{(m_1,t_1),\ldots,(m_n,t_n)\},\quad t_1<\cdots<t_n.
\ee
Then condition~(i) in the definition of an $[s,u]$-path in
Section~\ref{S:Gray} is equivalent to
\begin{itemize}
\item[(i)'] $m_k\big(\pi(t_k-)\big)\geq\pi(t_k)$ for all $k=1,\ldots,n$,
\item[(i)''] $t\mapsto\pi(t)$ is nondecreasing on
  $[s,t_1),\ldots,[t_{n-1},t_n),[t_n,u]$.
\end{itemize}
If $\pi$ would strictly increase at some time $t_{k-1}<t<t_k$ (with
$t_0:=s$ and $t_{n+1}:=u$), then we could make it smaller on $[t,t_k)$,
violating minimality (see (ii) in the definition of an $[s,u]$-path), so we
see that an $[s,u]$-path $\pi$ must actually satisfy
\be\label{between}
\mbox{$t\mapsto\pi(t)$ is constant on $[s,t_1),\ldots,[t_{n-1},t_n),[t_n,u]$.}
\ee
Condition~(i)' says that $\pi(t_k-)\in m_k^{-1}(\{\pi(t_k)\}^\up)$, 
which again by minimality implies that
\be\label{at}
\pi(t_k-)\in\big(m_k^{-1}\big(\{\pi(t_k)\}^\up\big)\big)_{\rm min}
\quad\Leftrightarrow\quad\pi(t_k-)\in m^{\dott}(\{\pi(t_k)\})
\qquad(k=1,\ldots,n),
\ee
where we have also used (\ref{mstar2}). By induction, decreasing $s$ wile
keeping $u$ fixed, it follows that $\zeta_{s,u}(B)\sub\Yb^\dott_{-u,-s}(B)$.

Conversely, one can verify that any cadlag function $\pi:[s,u]\to S$
satisfying (\ref{between}) and (\ref{at}) is an $[s,u]$-path.
This means that the function $[s,u]\ni t\mapsto\zeta_{t,u}$
is constant between the times $t_1,\ldots,t_n$ and satisfies
\be
\zeta_{t_k-,u}(B)
=\bigcup_{z\in\zeta_{t_k,u}(B)}\big(m_k^{-1}\big(\{z\}^\up\big)\big)_{\rm min}.
\ee
Again using (\ref{mstar2}), it follows that $[s,u]\ni t\mapsto\zeta_{t,u}$ is the
right-continuous modification of $[s,u]\ni t\mapsto\Yb^\dott_{-u,-t}$, proving
(\ref{Zaim}).
\end{Proof}

\section{Infinite product spaces}\label{S:infprod}

\subsection{Graphical representations}

In the present section we show how the additive and monotone systems dualities
from Sections~\ref{S:addual} and \ref{S:mondual} generalize to interacting
particle systems, which are processes with a state space of the form $T^\La$,
equipped with the product order, where $T$ is a finite partially ordered set
and $\La$ is a countably infinite set. Since we are interested in pathwise
dualities, we need a Poisson construction of such processes, which in this
context is usually called a \emph{graphical representation}. While the
generator construction of interacting particle systems is well-known (see
\cite{Lig85}), it is harder to find a good general reference for the graphical
construction. We base ourselves on the lecture notes \cite{Swa16}.

As in the finite case (see Section~\ref{S:randmap}), a graphical
representation of an interacting particle system is based on a collection
$\Gi$ of maps $m:T^\La\to T^\La$ together with rates $(r_m)_{m\in\Gi}$. As in
the finite case, one lets $\De$ be a Poisson point subset of $\Gi\times\R$
with local intensity $r_m\di t$, and defines $\De_{s,u}$ as before. Contrary
to the finite setting, however, $\De_{s,u}$ is typicaly an infinite set so one
cannot order the elements $(m,t)$ of $\De_{s,u}$ according to the time when
they occur and define a stochastic flow $(\Xb_{s,u})_{s\leq u}$ as in
(\ref{Phidef}). Instead, one defines
\be\label{finappr}
\Xb_{s,u}:=\lim_{\ti\De_n\up\De_{s,u}}\Xb^{\ti\De_n}_{s,u},
\ee
where for any finite $\ti\De\sub\De_{s,u}$
\be\ba{l}
\dis\Xb^{\ti\De}_{s,u}:=m_n\circ\cdots\circ m_1\\[5pt]
\dis\quad\mbox{with}\quad
\ti\De=\{(m_1,t_1),\ldots,(m_n,t_n)\},
\quad t_1<\cdots<t_n.
\ec
In order for $\Xb_{s,u}$ to be well-defined by (\ref{finappr}), one needs that
the limit exists and does not depend on the choice of the finite sets
$\ti\De_n\up\De_{s,u}$. For this, one needs to impose conditions on the maps
$m\in\Gi$ and their rates $r_m$.

For any function $m:T^\La\to T^\La$, we say that a point $j\in\La$
is \emph{$m$-relevant} for $i\in\La$ if
\be
\exists x,y\in T^\La\mbox{ s.t.\ }m(x)(i)\neq m(y)(i)
\mbox{ and }x(k)=y(k)\ \forall k\neq j,
\ee
i.e., changing the value of $x$ in $j$ may change the value of $m(x)$ in
$i$. We will use the notation
\be
\Ri_i(m):=\big\{j\in\La:j\mbox{ is $m$-relevant for }i\big\}.
\ee
A proof of the following lemma can be found in
\cite[Lemma~24]{SS15}. (Contrary to what the name ``$m$-relevant'' suggests,
condition~(ii) below may fail to be satisfied for discontinuous maps.)

\bl[Continuous maps]\label{L:mcont}
A map $m:T^\La\to T^\La$ is continuous with respect to the product topology on
$T^\La$ if and only if for each $i\in\La$, the following two conditions are
satisfied.
\begin{enumerate}
\item The set $\Ri_i(m)$ is finite.
\item If $x,y\in T^\La$ satisfy $x(j)=y(j)$ for all $j\in\Ri_i(m)$, then
  $m(x)(i)=m(y)(i)$.
\end{enumerate}
\el

For any map $m:T^\La\to T^\La$, let
\be
\Di(m):=\big\{i\in\La:\exists x\in T^\La\mbox{ s.t.\ }m(x)(i)\neq x(i)\big\}
\ee
denote the set of underlying space points whose values can possibly be changed
by $m$. We say that a map $m:T^\La\to T^\La$ is \emph{local} if $m$ is
continuous with respect to the product topology and $\Di(m)$ is finite.
We cite the following result from \cite[Thm~4.14]{Swa16}.

\bt[Well-defined graphical representation]\label{T:welldef}
Let $T$ be finite and $\La$ be countable. Let $\Gi$ be a countable collection
of local maps $m:T^\La\to T^\La$ and let $(r_m)_{m\in\Gi}$ be nonnegative
rates such that
\be\label{sum01}
\sup_{i\in\La}\sum\subb{m\in\Gi}{\Di(m)\ni i}r_m\big(|\Ri_i(m)|+1\big)<\infty.
\ee
Then almost surely for all $s\leq u$, the limit in (\ref{finappr}) exists in a
pointwise sense and does not depend on the choice of the finite sets
$\ti\De_n\up\De_{s,u}$.
\et

Under the condition (\ref{sum01}), it can moreover be shown (see
\cite[Thm~4.17]{Swa16}) that the stochastic flow $(\Xb_{s,u})_{s\leq u}$
defined in (\ref{finappr}) corresponds to a Feller process with generator of
the form (\ref{GPois}), defined using the usual generator construction of
interacting particle systems as in \cite[Thm~I.3.9]{Lig85}. In what follows,
we will in addition need the following lemma, which follows from
\cite[Lemma~4.12]{Swa16}.

\bl[Exponential bound]\label{L:expbd}
The maps $\Xb_{s,u}:T^\La\to T^\La$ are a.s.\ continuous for all $s\leq u$.
For deterministic $s\leq u$ and $i\in\La$, one has the exponential bound
\be
\E\big[\big|\Ri_i(\Xb_{s,u})\big|\big]\leq\ex{K(u-s)},
\ee
where $K$ is defined by
\be\label{summ}
K:=\sup_{i\in\La}\sum\subb{m\in\Gi}{\Di(m)\ni i}
r_m\big(|\Ri_i(m)|-1\big).
\ee
\el

\subsection{Monotone particle system duality}\label{S:monpart}

In the present section, we generalize Theorem~\ref{T:mondual} to monotonically
representable interacting particle systems with state space of the form
$S=T^\La$, where $T$ is a finite partially ordered set and $\La$ is
countable. For reasons that will become clear below, we assume that $T$ is
bounded from above with upper bound denoted by $1$. We equip $T^\La$ with the
product order $x\leq y$ iff $x(i)\leq y(i)$ for all $i\in\La$. Then $T^\La$ is
also bounded from above, with upper bound $1$ given by $1(i):=1$ $(i\in\La)$.
If $T'$ is a dual of $T$ as defined in Section~\ref{S:addual}, then ${T'}^\La$
is in a natural way a dual of $T^\La$, where we define the map $x\mapsto x'$
in a pointwise way as $x'(i):=(x(i))'$. Since $T$ is bounded from above, $T'$
is bounded from below with lower bound $0:=1'$. We also write $0$ to denote
the constant function $0(i):=0$ $(i\in\La)$ that is the lower bound of $T'$.

For any $x\in {T'}^\La$ and $B\in\Pc\big({T'}^\La\big)$, we let
\bc\label{suppx}
\dis{\rm supp}(x)&:=&\dis\{i\in\La:x(i)\neq 0\},\\[5pt]
\dis{\rm supp}(B)&:=&\dis\{i\in\La:x(i)\neq 0\mbox{ for some }x\in B\}
=\bigcup_{x\in B}{\rm supp}(x)
\ec
denote the ``support'' of $x$ and $B$, respectively. We let
\be\label{Sloc}
{T'}^\La_{\rm loc}:=\big\{x\in{T'}^\La:|{\rm supp}(x)|<\infty\}
\ee
denote the set of finitely supported $x\in {T'}^\La$, and write
$\Pc_{\rm fin}\big({T'}^\La_{\rm loc}\big)$ for the set of finite subsets of
${T'}^\La_{\rm loc}$. Equivalently,
\be
\Pc_{\rm fin}\big({T'}^\La_{\rm loc}\big)
=\big\{B\in\Pc\big({T'}^\La\big):|{\rm supp}(B)|<\infty\big\}.
\ee

As the state spaces for the processes $Y^\ast$ and $Y^\dgg$ from
Theorem~\ref{T:mondual}, we will choose
\be\label{Pcirc}
P_\ast:=\Pc_{\rm fin}({T'}^\La_{\rm loc})
\quand
P_\dgg:=\{B\in P_\ast:B=B_{\rm min}\}.
\ee
Note that these sets are countable, so the dual processes will be
continuous-time Markov chains. For any local map $m:T^\La\to T^\La$ and $B\in
P_\ast$, we define $m^\dgg(B)$ and $m^\ast(B)$ as in (\ref{mast}), and we
define a duality function $\phi:T^\La\times P_\ast\to\{0,1\}$ as in
(\ref{phidef}).

\bl[Duals of monotone local maps]\label{L:infdumaps}
Let $T$ be a finite partially ordered set that is bounded from above, let
$\La$ be countable, and let $m:T^\La\to T^\La$ be a local map that is
monotone. Then the maps $m^\ast$ and $m^\dgg$ defined in (\ref{mast}) map the
space $P_\ast$ into itself, and are dual to $m$ with respect to the duality
function $\phi$ from (\ref{phidef}). Moreover, (\ref{dagprop}) and
(\ref{cuprop}) hold for all $B,C\in P_\ast$. In particular, by
(\ref{dagprop}), $m^\dgg$ maps $P_\ast$ into $P_\dgg$.
\el
\begin{Proof}
Let us say that a set $A\sub\Pc({T'}^{\Lambda})$ is \emph{locally
defined} if there exists a finite set $\Ga\sub\La$ and a set $C\sub{T'}^\Ga$
such that
\be
A=C\times{T'}^{\La\beh\Ga}.
\ee
For a set of this form, by Lemma~\ref{L:minmax}
\be
\big(C\times{T'}^{\La\beh\Ga}\big)_{\rm min}=C_{\rm min}\times\{0\}
\quand
\big(\big(C\times{T'}^{\La\beh\Ga}\big)_{\rm min}\big)^\up
=C_{\rm min}^\up\times\{0\}^\up\supset C\times{T'}^{\La\beh\Ga},
\ee
where $0$ here denotes the minimal element of ${T'}^{\La\beh\Ga}$.
In particular, this shows that for any locally defined increasing set $A$,
\be
A_{\rm min}\in P_\ast
\quand
(A_{\rm min})^\up=A.
\ee
Also, obviously, $B\in P_\ast$ implies that $B^\up$ is locally defined.
The proof of Lemma~\ref{L:dumaps} now carries over without a change, where we
use that since $m$ is a local map, the map $n:{T'}^{\La} \to {T'}^{\La}$
defined in (\ref{nacm}) has the property that
\be
A\mbox{ locally defined implies }n^{-1}(A)\mbox{ locally defined}
\qquad\big(A\in\Pc({T'}^{\La})\big).
\ee
\end{Proof}

Let $\Gi$ be a collection of monotone local maps $m:T^\La\to T^\La$, let
$(r_m)_{m\in\Gi}$ be nonnegative rates, and let $\De$ be a Poisson point
subset of $\Gi\times\R$ with local intensity $r_m\di t$. As in
(\ref{hatDe}), we define
\be\label{hatDe2}
\De^\ast:=\{(m^\ast,-t):(m,t)\in\De\}
\quand
\De^\dgg:=\{(m^\dgg,-t):(m,t)\in\De\}.
\ee
For $\both=\ast$ or $\dgg$, we wish to define
\be\label{Yfinappr}
\Yb^\both_{s,u}:=\lim_{\ti\De_n\up\De^\both_{s,u}}\Yb^{\both,\,\ti\De_n}_{s,u},
\ee
where for any finite $\ti\De\sub\De^\both_{s,u}$
\be\ba{l}
\dis\Yb^{\both,\,\ti\De}_{s,u}:=m^\both_n\circ\cdots\circ m^\both_1\\[5pt]
\dis\quad\mbox{with}\quad
\ti\De=\{(m^\both_1,t_1),\ldots,(m^\both_n,t_n)\},
\quad t_1<\cdots<t_n.
\ec
It turns out that this is all right under the same technical condition
(\ref{sum01}) that we used to prove the the forward process is well-defined.

\bp[Well-definedness of the dual processes]\label{P:Zdef}
Under the condition (\ref{sum01}) almost surely for all $s\leq u$, the limit
in (\ref{Yfinappr}) exists in a pointwise sense and does not depend on the
choice of the finite sets $\ti\De_n\up\De^\both_{s,u}$. Moreover,
\be\label{Esup}
\E\big[\big|{\rm supp}\big(\Yb^\both_{s,u}(B)\big)\big|\big]
\leq|{\rm supp}(B)\big|\ex{K(u-s)}\qquad(s\leq u),
\ee
where $K$ is the constant in (\ref{summ}).
\ep

The following lemma prepares for the proof of Proposition~\ref{P:Zdef}.

\bl[Support and dual maps]
Let $\both=\ast$ or $\dgg$ and let $m:T^\La\to T^\La$ be a local map.
Then, for any $B\in P_\both$,
\be\label{inzet}
{\rm supp}\big(m^\both(B)\big)\sub
\bigcup_{i\in{\rm supp}(B)}\Ri_i(m)
\ee
and
\be\label{supmatters}
{\rm supp}(B)\cap\Di(m)=\emptyset
\quad\mbox{implies}\quad
m^\both(B)=B.
\ee
\el
\begin{Proof}
By (\ref{dagprop}), which by Lemma~\ref{L:infdumaps}
holds for all $B\in P_\ast$, we have $m^\both(B)\sub m^\ast(B)$ $(B\in
P_\ast)$, so it suffices to prove (\ref{inzet}) for $m^\ast$ only. 
By (\ref{cuprop}), which by Lemma~\ref{L:infdumaps} holds for all $B\in
P_\ast$, both the left- and right-hand side of (\ref{inzet}) are additive as a
function of $B$, so it suffices to prove the statement for one-point sets of
the form $B=\{x\}$ with $x\in{T'}^\La_{\rm loc}$. In this case,
\be
m^\ast\big(\{x\}\big)=n^{-1}\big(\{x\}^\up\big)_{\rm min},
\ee
where $n$ is defined in terms of $m$ as in (\ref{nacm}). In other words, this
says that
\be\label{mz}
m^\ast\big(\{x\}\big)=\big\{y\in{T'}^\La:m(y')'\geq x\big\}_{\rm min}
=\big\{y\in{T'}^\La:m(y')\leq x'\big\}_{\rm min}.
\ee
For each $k\in{\rm supp}\big(m^\ast\big(\{x\}\big)\big)$ there exists a
$z\in m^\ast\big(\{x\}\big)$ such that $z(k)\neq 0$. Define $y$ by $y(k):=0$
and $y(j)=z(j)$ for all $j\neq k$. Then, by (\ref{mz}), $m(z')\leq x'$ but
$m(y')\not\leq x'$ by the minimality of $z$. It follows that
$k\in\Ri_i(m)$ for some $i$ such that $x'(i)\neq 1$, proving that
\be
{\rm supp}\big(m^\ast(\{x\})\big)\sub
\bigcup_{i\in{\rm supp}(\{x\})}\Ri_i(m).
\ee

To prove also (\ref{supmatters}), we observe that ${\rm
  supp}(B)\cap\Di(m)=\emptyset$ implies $m^{-1}({B'}^\down)={B'}^\down$ and
hence, by (\ref{mast})
\be
m^\dgg(B)':=({B'}^\down)_{\rm max}
\quad\mbox{which implies}\quad
m^\dgg(B)=(B^\up)_{\rm min}=B_{\rm min}=B\qquad(B\in P_\dgg)
\ee
and similarly
\be
m^\ast(B):=\bigcup_{x\in B}(\{x\}^\up)_{\rm min}=\bigcup_{x\in B}\{x\}=B
\qquad(B\in P_\ast).
\ee
\end{Proof}

\begin{Proof}[of Proposition~\ref{P:Zdef}]
Fix $B\in P_\both$ and $s\leq u$. Set $t_0:=s$, $\Ga_0:={\rm supp}(B)$,
inductively define times $t_k$ and sets $\Ga_k$ $(k\geq 1)$ by
\be
t_k:=\inf\big\{t>t_{k-1}:\exists\ (m^\both,t)\in\De^\both_{s,u}\mbox{ s.t.\ }
\Ga_{k-1}\cap\Di(m^\both)\neq\emptyset\big\},
\ee
and let $m^\both_k$ denote the unique element of $\Gi$ such that
$(m^\both_k,t_k)\in\De^\both_{s,u}$. Then, exactly as in the proof of
\cite[Lemma~4.13]{Swa16}, using (\ref{inzet}), one can see that (\ref{sum01})
implies that this inductive process ends after a finite number of steps, i.e.,
there exists some $n$ such that $\Ga_n\cap\Di(m^\both)=\emptyset$ for all
$t>t_n$ and $(m^\both,t)\in\De^\both_{s,u}$.

By (\ref{supmatters}), as soon as the finite sets $\ti\De_n$ from
(\ref{Yfinappr}) contain the elements
$(m^\both_1,t_1),\ldots,$ $(m^\both_n,t_n)$, adding more elements to $\ti\De_n$
has no effect on $\Yb^{\both,\,\ti\De_n}_{s,u}(B)$ and hence the limit in
(\ref{Yfinappr}) exists independent of the choice of the sequence $\ti\De_n$.
Now (\ref{inzet}) implies that
\be
{\rm supp}\big(\Yb^\both_{s,u}(B)\big)\sub
\bigcup_{i\in{\rm supp}(B)}\Ri_i(\Xb_{-u,-s})
\ee
and hence (\ref{Esup}) follows from Lemma~\ref{L:expbd}.
\end{Proof}

Using Proposition~\ref{P:Zdef}, it is straightforward to check that
$(\Yb^\both_{s,u})_{s\leq u}$ is a stochastic flow with independent
increments and that if $B_0$ is a $P_\both$-valued random variable,
independent of the Poisson set $\De^\both$, then
\be
Y^\both_t:=\Yb^\both_{0,t}(B_0)\qquad(t\geq 0)
\ee
defines an (obviously nonexplosive) Markov process $(Y^\both_t)_{t\geq
  0}$ with countable state space $P_\both$ and generator $H_\both$ as in
(\ref{HPoisY}) for $\both=\ast$ or $\dgg$.

\bp[Pathwise duality for monotone particle systems]\label{P:monpart}
Let $T$ be a finite partially ordered set that is bounded from above, let
$\La$ be countable, and let $X$ be an interacting particle system whose
generator has a random mapping representation of the form (\ref{GPois}),
where all elements of $\Gi$ are monotone local maps and the rates satisfy the
summability condition (\ref{sum01}). Let $\De$ be a graphical representation
for $X$ and define a stochastic flow $(\Xb_{s,u})_{s\leq u}$ as in
(\ref{finappr}). For $\both=\ast$ or $\dgg$, define a Poisson set $\De^\both$
as in (\ref{hatDe2}) and use this to define random maps
$(\Yb^\both_{s,u})_{s\leq u}$ on the space $P_\both$ from (\ref{Pcirc})
as in (\ref{Yfinappr}). Then $(\Yb^\both_{s,u})_{s\leq u}$ is a stochastic
flow that is dual to $(\Xb_{s,u})_{s\leq u}$ with respect to the 
duality function $\phi$ in (\ref{phidef}), in the sense defined in
Section~\ref{S:pathdual}. In particular, for each $s<u$, $x\in T^\La$, and
$B\in P_\both$, the function
\be
[s,u]\ni t\mapsto1_{\txt\{\Xb_{s,t-}(x)\leq z\mbox{ for some }z\in\Yb^\both_{-u,-t}(B)\}}
\ee
is a.s.\ constant.
\ep
\begin{Proof}
Using Lemma~\ref{L:infdumaps}, this follows just as in the proof of
Proposition~\ref{P:pathdual}.
\end{Proof}

\subsection{Additive particle systems}\label{S:adpart}

In the present subsection, specializing from the set-up of the previous
section, we look at interacting particle systems that are defined by additive
local maps and whose state space is of the form $T^\La$ with $T$ a finite
lattice. As before, we equip $T^\La$ with the product order; then
  $T^\La$ is also a lattice, where $(x\vee y)(i)=x(i)\vee y(i)$
and $(x\wedge y)(i)=x(i)\wedge y(i)$. We let $T'$
denote a dual of $T$ so that ${T'}^\La$ is in a natural way a dual of $T^\La$.
We recall from (\ref{tetSS}) that
\be\label{tetSS2}
\li x,y\re:=1_{\txt\{x\leq y'\}}=1_{\txt\{y\leq x'\}}
\qquad(x\in T^\La,\ y\in {T'}^\La).
\ee
The following lemma generalizes Lemmas~\ref{L:addual} and
\ref{L:astac} to infinite product spaces.

\bl[Additive local maps]\label{L:locdu}
Let $T$ be a finite lattice, let $T'$ be its dual, and let $\La$ be a
countable set. Then, for each additive local map $m:T^\La\to T^\La$ there
exists a unique additive local map $m':{T'}^\La\to{T'}^\La$ such that
\be\label{dumap3}
\li m(x),y\re=\li x,m'(y)\re
\qquad(x\in T^\La,\ y\in {T'}^\La),
\ee
where $\li\,\cdot\,,\,\cdot\,\re$ is as defined in (\ref{tetSS}).
For $\both=\ast$ or $\dgg$, let $m^\both:P_\both\to P_\both$ be defined as in
(\ref{mast}). Then
\be\label{astac2}
m^\ast(B)=\{m'(y):y\in B\}
\quand
m^\dgg(B)=m^\ast(B)_{\rm min}
\ee
for all $B\in P_\ast$ resp.\ $B\in P_\dgg$.
\el
\begin{Proof}
Let $\ti\La:=\Di(m)\cup\bigcup_{i\in\Di(m)}\Ri_i(m)$. Then $\ti\La$ is a finite
set by the definition of a local map and Lemma~\ref{L:mcont}. Moreover, by the
same lemma, there exists an additive map $\ti m:T^{\ti\La}\to T^{\ti\La}$ such
that
\be
m(x)(i)=\left\{\ba{ll}
\dis\ti m\big((x(j))_{j\in\ti\La}\big)(i)\quad&\mbox{if }i\in\ti\La,\\[5pt]
x(i)\quad&\mbox{otherwise.}\ea\right.
\ee
By Lemma~\ref{L:addual}, $\ti m$ has a unique dual map $\ti m'$, which
is also additive. Lifting this map to the larger space yields a local map $m'$
as in (\ref{dumap3}). Since knowing $\li x,m'(y)\re$ for all $x\in T^\La$
determines $m'(y)$ uniquely, such a map is unique.

The proof of formula (\ref{astac2}) is the same as in the finite case
(Lemma~\ref{L:astac}).
\end{Proof}

Let $\Gi$ be a collection of additive local maps $m:T^\La\to
T^\La$ and let $(r_m)_{m\in\Gi}$ be nonnegative rates.
Apart from the summability condition (\ref{sum01}) we will sometimes also need
the dual condition
\be\label{dusum01}
\sup_{i\in\La}\sum\subb{m\in\Gi}{\Di(m')\ni i}r_m\big(|\Ri_i(m')|+1\big)<\infty.
\ee
For $x\in T^\La$, we define ${\rm supp}(x)$ as in (\ref{suppx}) and we define
$T^\La_{\rm loc}$ as in (\ref{Sloc}) (but with $T'$ replaced by $T$).
We make the following observation.

\bl[Finite initial states]\label{L:fininf}
Assuming the dual condition (\ref{dusum01}) but not necessarily (\ref{sum01}),
almost surely for all $s\leq u$, the limit in (\ref{finappr}) exists in a
pointwise sense on $T^\La_{\rm loc}$, does not depend on the choice of the
finite sets $\ti\De_n\up\De_{s,u}$, and defines a random map
$\Xb_{s,u}:T^\La_{\rm loc}\to T^\La_{\rm loc}$.
\el
\begin{Proof}
We may equivalently show that if the rates $(r_m)_{m\in\Gi}$ satisfy the
summability condition (\ref{sum01}), then the dual graphical representation
$\De':=\{(m',-t):(m,t)\in\De\}$ unambiguously defines random maps 
$\Yb_{s,u}:{T'}^\La_{\rm loc}\to {T'}^\La_{\rm loc}$.
By formula (\ref{astac2}), the map $m^\ast$ maps the space of all
singleton-sets $\{y\}$ with $y\in{T'}^\La_{\rm loc}$ into itself, and
\be
m^\ast\big(\{y\}\big)=\{m'(y)\}\qquad(y\in{T'}^\La_{\rm loc}).
\ee
In view of this, the statement follows from Proposition~\ref{P:Zdef}.
\end{Proof}

The following theorem generalizes Theorem~\ref{T:addual} to interacting
particle systems.

\bt[Duality for additive particle systems]
Let $T$ be a finite lattice and let $\La$ be a countable set. Let $\Gi$ be a
collection of additive local maps $m:T^\La\to T^\La$ and let $(r_m)_{m\in\Gi}$
be nonnegative rates satisfying the summability condition (\ref{sum01}).  Let
$\De$ be a Poisson point subset of $\Gi\times\R$ with local intensity $r_m\di
t$ and let $\De':=\{(m',-t):(m,t)\in\De\}$. Let $(\Xb_{s,t})_{s\leq t}$ and
$(\Yb_{s,t})_{s\leq t}$ be the stochastic flows on $T^\La$ and ${T'}^\La_{\rm
  loc}$ defined from $\De$ and $\De'$ by grace of Theorem~\ref{T:welldef} and
Lemma~\ref{L:fininf}, respectively. Then $(\Xb_{s,t})_{s\leq t}$ and
$(\Yb_{s,t})_{s\leq t}$ are dual with respect to the  duality function
$\psi(x,y):=\li x,y\re$, in the sense defined in Section~\ref{S:pathdual}. In
particular, for each $x\in T^\La$, $y\in {T'}^\La_{\rm loc}$, and $s\leq u$,
the
function
\be\label{mappsi2}
[s,u]\mapsto\li\Xb_{s,t-}(x),\Yb_{-u,-t}(y)\re
\ee
is a.s.\ constant. If moreover the dual summability condition (\ref{dusum01})
is satisfied, then the same statements hold for the stochastic flow
$(\Yb_{s,t})_{s\leq t}$ on the larger state space ${T'}^\La$
(instead of ${T'}^\La_{\rm loc}$).
\et
\begin{Proof}
As explained in the proof of Lemma~\ref{L:fininf}, by (\ref{astac2}), for
$y\in{T'}^\La_{\rm loc}$, the dual map $\Yb_{s,t}$ is a special case of
the dual map $\Yb^\ast_{s,t}$ from the previous subsection, so the
first claim follows from Proposition~\ref{P:monpart}. If moreover the dual
  summability condition (\ref{dusum01}) is satisfied, then by
Theorem~\ref{T:welldef}, $\Yb_{s,t}$ is also well-defined on
  ${T'}^\La$. For any $y\in{T'}^\La$, we can find $y_n\in{T'}^\La_{\rm loc}$
that increase to $y$. It is not hard to see that this implies that
$\Yb_{s,t}(y_n)$ increases to $\Yb_{s,t}(y)$ and hence, for all $t\in[s,u]$,
\be
\li\Xb_{s,t-}(x),\Yb_{-u,-t}(y_n)\re
=1_{\txt\{\Yb_{-u,-t}(y_n)\leq\Xb'_{s,t-}(x)\}}\down
\li\Xb_{s,t-}(x),\Yb_{-u,-t}(y)\re
\ee
so the limit is a.s.\ constant as a function of $t$.
\end{Proof}

\appendix

\section{A bit of lattice theory}\label{S:bit}

In this appendix we collect some elementary properties of partially ordered
sets and lattices that are used in the paper.

Let $S$ be a any set. Recall that a relation $\leq$ on $S$ is called a
\emph{partial order} if it satisfies the axioms
\begin{enumerate}
\item $x\leq x$ $(x\in S)$.
\item $x\leq y$ and $y\leq x$ implies $x=y$ $(x,y\in S)$.
\item $x\leq y\leq z$ implies $x\leq z$ $(x,y\in S)$.
\end{enumerate}
A \emph{partially ordered set} (also called \emph{poset}) is a set with a
partial order defined on it. A \emph{total order} is a partial order such that
moreover
\begin{enumerate}\addtocounter{enumi}{3}
\item $x\leq y$ or $y\leq x$ (or both) $(x,y\in S)$.
\end{enumerate}

Increasing and decreasing sets and the notation $A^\up$ and $A^\down$ have
already been defined in Section~\ref{S:subspace}.

By definition, a \emph{minimal} element of a set $A\sub S$ is an $x\in A$ such
that $\{y\in A:y\leq x\}=\{x\}$. \emph{Maximal} elements are minimal elements
for the reversed order. The following simple observation is well-known.

\bl[Maximal elements]\label{L:minmax}
Let $S$ be a partially ordered set and let $A\sub S$ be finite.
Then, for every $x\in A$ there exists a maximal element $y\in A$ such that
$y\geq x$.
\el
\begin{Proof}
Either $x$ is a maximal element or there exists some $x'\in A$, $x'\neq x$,
such that $x'\geq x$. Continuing this process, using the finiteness of $A$, we
arrive after a finite number of steps at a maximal element.
\end{Proof}

In particular, Lemma~\ref{L:minmax} shows that every nonempty finite $A\sub S$
has a maximal element. Applying this to the reversed order, we see that $A$
also contains a minimal element.

Let $S$ be a partially ordered set and $A\sub S$. An \emph{upper bound} of $A$
is an element $x\in S$ such that $y\leq x$ for all $y\in A$. Note that in
general $x$ does not have to be an element of $A$. If there exists an element
$x\in A$ that is an upper bound of $A$, then such an element is necessarily
unique. (Indeed, imagine that $x'\in A$ is also an upper bound for $A$, then
$x'\leq x$ and $x\leq x'$.)

If $A$ is decreasing and there exists an element $x\in A$ that is an upper
bound of $A$, then $A=\{x\}^\down$. Indeed, the fact that $x$ is an upper
bound means that $A\sub\{x\}^\down$ while by the facts that $x\in A$ and $A$
is decreasing, we have $A\supset\{x\}^\down$. A \emph{least upper bound} of
$A$ is an element $x\in S$ that is an upper bound of $A$ and that satisfies
$x\leq x'$ for any (other) upper bound $x'$ of $A$. Note that $\bigcap_{x\in
  A}\{x\}^\up$ is the set of all upper bounds of $A$, which is an increasing
set. Thus, a least upper bound of $A$ is an element $x\in\bigcap_{x'\in
  A}\{x'\}^\up$ that is a lower bound of $\bigcap_{x'\in A}\{x'\}^\up$. By our
earlier remarks, such an element is unique. Also, since $\bigcap_{x'\in
  A}\{x'\}^\up$ is an increasing set, if there exists an element
$x\in\bigcap_{x'\in A}\{x'\}^\up$ that is a lower bound of $\bigcap_{x'\in
  A}\{x'\}^\up$, then $\bigcap_{x'\in A}\{x'\}^\up=\{x\}^\up$.

A partially ordered set $S$ is a join-semilattice if and only if for each
$y,z\in S$, the set $\{y,z\}$ has a least upper bound. Equivalently, this says
that there exists an element $x\in\{y\}^\up\cap\{z\}^\up$ that is a lower
bound of $\{y\}^\up\cap\{z\}^\up$. By our earlier remarks, such an element is
unique and one must in fact have $\{y\}^\up\cap\{z\}^\up=\{x\}^\up$. We denote
this unique element by $x=:y\vee z$ and call it the \emph{supremum} or
\emph{join} of $x$ and $y$. Note that this (more traditional) definition of
suprema and join-semilattices is equivalent to the definition in
(\ref{supdef}). Infima and meet-semilattices are defined in the same way, but
with respect to the reversed order.

Each finite join-semilattice $S$ is clearly bounded from above, since
$\bigvee_{x\in S}x$ is an upper bound; similarly finite meet-semilattices are
bounded from below.

As already mentioned in Section~\ref{S:subspace}, if $S$ is a partially
ordered set, then a nonempty increasing set $A\sub S$ such that for every
$x,y\in A$ there exists a $z\in A$ with $z\leq x,y$ is called a \emph{filter}
and a nonempty decreasing set $A$ such that for every $x,y\in A$ there exists
a $z\in A$ with $x,y\leq z$ is called an \emph{ideal}.

The following lemma characterizes ideals in join-semilattices. An analogue
statement holds for filters. In a lattice-theoretic setting, this lemma is
often taken as the definition of an ideal.

\bl[Ideals in join-semilattices]\label{L:latfil}
Let $S$ be a join-semilattice and let $A\sub S$. Then $A$ is an ideal if and
only if $A$ is nonempty, decreasing, and $x,y\in A$ imply $x\vee y$ in $A$.
\el
\begin{Proof}
Let $A$ be nonempty and decreasing. If $x,y\in A$ imply $x\vee y$ in $A$, then
for every $x,y\in A$ there is an element of $A$, namely $x\vee y$, such that
$x\vee y\geq x,y$, showing that $A$ is an ideal. Conversely, if $A$ is an
ideal, then for every $x,y\in A$ there is an element $z\in A$ such that $z\geq
x,y$ and hence $z\geq x\vee y$, which by the fact that $A$ is decreasing
proves that $x\vee y\in A$.
\end{Proof}

As already mentioned in Section~\ref{S:subspace}, a \emph{principal filter} is
a filter that contains a minimal element and a \emph{principal ideal} is an
ideal that contains a maximal element.  We state the following facts for
ideals only; applying them to the reversed order yields analogue statements
for filters.

\bl[Principal ideals]\label{L:filter}
Let $S$ be a partially ordered set and let $A\sub S$. Then:
\begin{itemize}
\item[{\rm 1.}] If $A$ is an ideal, then $A$ contains at most one maximal
  element.
\item[{\rm 2.}] $A$ is a principal ideal if and only if there exists some
  $z\in S$ such that $A=\{z\}^\down$.
\end{itemize}
If $A$ is moreover finite, then the following conditions are equivalent:
\begin{enumerate}
\item $A$ is a principal ideal.
\item $A$ is an ideal.
\item $A$ is decreasing and contains a unique maximal element.
\end{enumerate}
\el
\begin{Proof}
1.\ If $A$ is an ideal and $x,y\in A$ are maximal, then there exists some
$z\in A$ such that $z\geq x,y$ and hence $x=z=y$ by the definition of
maximality.

2.\ Let $A$ be a principal ideal with maximal element $z$. Then
$\{z\}^\down\sub A$ by the fact that $A$ is decreasing. The inclusion
$\{z\}^\down\supset A$ follows by observing that if $x\in
A$, then by the definition of an ideal there exists some $y\in A$ with $y\leq
x,z$, and hence $y=z$ by the maximality of $z$, so $x\leq y=z$.
Conversely, it is straightforward to check that each set of the form
$\{z\}^\down$ is a principal ideal.

%Indeed, $\{z\}^\down$ is clearly nonempty and decreasing, while for every
%$x,y\in\{z\}^\down$ there exists an element of $\{z\}^\down$, namely $z$, such
%that $z\geq x,y$, proving that $\{z\}^\down$ is an ideal. Since $z$ is the
%maximal element of $\{z\}^\down$, this is a principal ideal.

%(Indeed, if $x\in\{z\}^\down$ satisfies $x\geq z$, then
%clearly $x=z$ proving that $z$ is a maximal element of $\{z\}^\up$.)

The implications (i)$\volgt$(ii) and (i)$\volgt$(iii) are trivial. Assume that
$A\sub S$ is finite. Then (ii)$\volgt$(i) by Lemma~\ref{L:minmax}. We claim
that moreover (iii)$\volgt$(i). Let $z$ be the unique maximal element of $A$. By
2, it suffices to prove that $A=\{z\}^\down$. Since $A$ is decreasing, we have
$\{z\}^\down\sub A$. To prove the other inclusion, we observe that for any $x\in
A$, by Lemma~\ref{L:minmax}, there exists a maximal element $y\in A$ such that
$y\geq x$. By assumption, $z$ is the unique maximal element of $A$, so $y=z$
and hence $x\leq z$, proving that $x\in\{z\}^\down$.
\end{Proof}

\newpage
\section{List of notation}\label{S:not}
\begin{itemize}
\item $\sub$ used in the meaning of $\subseteq$ throughout the paper.
\item $0,1$ general notation for the \emph{lower} and \emph{upper bound} of a
  partially ordered set (as defined below (\ref{supdef})).
\item $\li x,y\re:=1_{\txt\{x\leq y'\}}, x\in S,\ y\in S'$ (see (\ref{tetSS})).
\item $\leadsto$ indicates the existence of an \emph{open path} in a
  percolation substructure (see Section~\ref{S:percol}).
\item $\vee,\wedge$ \emph{supremum, infimum} (see (\ref{supdef})).
\item $A^\up:=\{x\in S:x\geq y\mbox{ for some }y\in A\}\supset A$. 
\item $A^\down:=\{x\in S:x\leq y\mbox{ for some }y\in A\}\supset A$. 
\item $A_{\rm max},A_{\rm min}$ \emph{set of maximal}, resp.\ \emph{minimal
  elements of A} (see (\ref{Amax})).
\item$\Fi(S,T)$ \emph{space of all functions $f:S\to T$}.
\item $\Fi_{\rm mon}(S,T)$ \emph{space of all monotone maps $m:S\to T$} (see (\ref{mondef})).
\item $\Gi$ set of maps used in the random mapping representation of a
  generator (see Section~\ref{S:randmap}).
\item $G, H$ \emph{Markov process generator and dual generator} (see
  (\ref{GPois}) and (\ref{HPois})).
\item $H_\dgg, H_\ast, H_\circ, H_\dott$ \emph{dual generators} (see
  (\ref{HPoisY})).
\item $\Ki(S,T)$ \emph{space of all probability kernels from $S$
to $T$}. %see (\ref{randmap})
\item $P_t$ \emph{transition probabilities} (see Section~\ref{S:pathdual}).
\item $\Pc(S):=\{A:A\sub S\}$ \emph{set of all subsets of $S$}.
\item $\dis\Pc_{\rm inc}(S),\dis\Pc_{\rm dec}(S)$ \emph{increasing, decreasing
  subsets of S} (see (\ref{!dec})).
\item $\dis\Pc_{\rm !inc}(S),\dis\Pc_{\rm !dec}(S)$ \emph{principal filters,
  ideals of S} (see (\ref{!dec})).
\item $S,T$ \emph{state spaces}, mostly endowed with a partial order $\leq.$
\item $S'$ \emph{dual of the partially ordered space S} (see (\ref{duspa})).
\item $X_t,Y_t$ \emph{Markov process} and \emph{dual Markov process} (see
  Section~\ref{S:pathdual}).
\item $\Xb_{s,t},\Yb_{s,t}$ \emph{stochastic flows} associated with Markov
  processes $X,Y$ (see Section~\ref{S:pathdual}).
\item $m:S\to S$ and $\hat m:T\to T$ \emph{map and dual map} with respect to
  $\psi$ (see (\ref{mdual})).
\item $m'$ \emph{dual map} associated with an additive map
  $m$ (see Lemma~\ref{L:addual}).
\item $m^\dgg, m^\ast, m^\circ, m^\dott$ \emph{dual maps} associated with a
  monotone map $m$ (see (\ref{mast}) and (\ref{mstar})).
\item $x^{\rm c}$ \emph{complement} of a set $x$.
\item $x',y'$ \emph{elements of the dual space} $S'$ associated with $x,y\in
  S$ (see (\ref{duspa})).
\item $\De,\De_{s,t}$ \emph{Poisson point sets} used in the \emph{graphical
  representation} of a Markov process (see Section~\ref{S:randmap}).
\item $\Lambda$ a \emph{finite set} (in
  Sections~\ref{S:intro}--\ref{S:proofs}) or a
  \emph{countable set} (in Section~\ref{S:infprod}).
\item $\phi, \ti\phi:$ \emph{duality functions} (see (\ref{phidef}) and
  (\ref{altphi})).
\item $\psi:S\times T\to\R$ \emph{duality function} (see (\ref{dual}) and
  (\ref{dual2})).
%\item $\ov m, \ov G$  
%\item $M, M'$
%\item $\gamma$
%lattice
%join-semilattice
%meet-semilattice
%bounded from below/above
%monotone map
%additive map
%supremum/join
%infimum/meet
\end{itemize}

\noindent
{\large {\bf Acknowledgements}}\\
%\section{Acknowledgements}
We thank Fero Mat\'u\v{s} and June Huh for useful discussions, and
L\'aszl\'o Csirmaz for help with the proof of Lemma~\ref{L:nondis}. We also
thank two anonymous referees for a careful reading of the manuscript and 
some helpful suggestions.



\end{document}